\documentclass[11pt,reqno]{amsart}
\usepackage[usenames]{color}
\usepackage{amsmath, verbatim,graphicx,epstopdf,enumerate}
\newcommand{\lb}{\left(}

\newcommand{\rb}{\right)}

\newcommand{\Lc}{\mathcal{L}}

\newcommand{\Beq}{\begin{equation}}
	\newcommand{\Eeq}{\end{equation}}
\newcommand{\beq}{\begin{equation*}}
	\newcommand{\eeq}{\end{equation*}}
\newcommand{\bal}{\begin{align}}
	\newcommand{\eal}{\end{align}}

\newcommand{\bp}{\begin{prob}}
	\newcommand{\ep}{\end{prob}}
\newcommand{\bpr}{\begin{proof}}
	\newcommand{\epr}{\end{proof}}

\newcommand{\bel}[1]{\begin{equation}\label{#1}}
	\newcommand{\ee}{\end{equation}}

\newcommand{\rr}{\mathbb{R}}

\newtheorem{theorem}{Theorem}[section]

\newtheorem{lemma}[theorem]{Lemma}
\newtheorem{prop}[theorem]{Proposition}

\theoremstyle{definition}
\newtheorem{definition}[theorem]{Definition}
\newtheorem{remark}[theorem]{Remark}

\usepackage[notref,notcite]{}
\usepackage[pagewise]{lineno}

\numberwithin{equation}{section}

\usepackage{amsthm}

\usepackage{amsfonts}
\newcommand{\vph}{\varphi}

\title[Partial data inverse problem]{Inverse problem for a time-dependent Convection-diffusion equation in admissible geometries}
\author[Mishra, Purohit and Vashisth]{Rohit Kumar Mishra$^\dagger$, Anamika Purohit$^\ddagger$ and Manmohan Vashisth$^{\ast}$}

\address{$^\dagger$ Department of Mathematics, Indian Institute of Technology, Gandhinagar, Gujarat - 382355, India.
\newline
\indent E-mail:{\tt \ rohit.m@iitgn.ac.in, \ rohittifr2011@gmail.com}}
\address{$^\ddagger$ Department of Mathematics, Indian Institute of Technology, Gandhinagar, Gujarat - 382355, India.
\newline
\indent E-mail:{\tt \ anamika.purohit@iitgn.ac.in }}
\address{$^\ast$Department of Mathematics,  Indian Institute of Technology Ropar, Rupnagar, Punjab - 140001, India.
\newline
\indent E-mail:{\tt\  manmohanvashisth@iitrpr.ac.in,\  manmohanvashisth@gmail.com}}

\begin{document}
\begin{abstract}
We consider a partial data inverse problem for a time-dependent convection-diffusion equation on an admissible manifold. We prove that the time-dependent convection term and time-dependent density can be recovered uniquely modulo a known gauge invariance. There have been several works on inverse problems related to the steady state convection-diffusion operator in Euclidean as well as in Riemannian geometry settings; however, inverse problems related to time-dependent convection-diffusion equation on a manifold are not studied in the prior works, which is the main aim of this paper.  In fact, to the best of our knowledge, the problem studied here is the first work related to a partial data inverse problem for recovering both first and zeroth-order time-dependent perturbations of evolution equations in the Riemannian geometry setting.
\end{abstract}
\maketitle


\section{Introduction and statement of main result} The paper deals with a partial data inverse problem related to a convection-diffusion equation on $M_{T}:=(0, T)\times M$ where $0<T<\infty$ and $(M,g)$ is a smooth $n-$dimensional ($n\geq 2$)  Riemannian manifold having smooth boundary $\partial M$. We denote by $\Sigma:=(0,T)\times \partial M$ as the lateral boundary of $M_{T}$ and $\displaystyle \partial M_{T}:=\Sigma\cup (\{0\}\times M)\cup (\{T\}\times M)$ the topological boundary of $M_T$. We also denote by ${T}M$ and ${T}^{*}M$ the tangent and cotangent bundle {of} $M$. For a convection term $A\in W^{1,\infty}\left(M_{T}; T^{*}M\right)$ given by $ \displaystyle A(t,x):=\sum_{j=1}^{n}A_{j}(t,x)dx^{j}$ in local coordinates $x_1,x_2,\dots, x_n$ of manifold $M$ and density $q\in L^{\infty}(M_T)$, the initial boundary value problem (IBVP)  for the convection-diffusion equation on $M_{T}$ is modeled by the following IBVP for second order linear parabolic partial differential equation (PDE)
\begin{align}\label{eq: Main equation}
\begin{aligned}
\begin{cases}
\Big[\partial_{t}-\sum\limits_{j,k=1}^{n}\frac{1}{\sqrt{\lvert g\rvert}}\Big(\partial_{x_{j}}+ A_{j}\Big)\left(\sqrt{\lvert g\rvert} g^{jk}\left (\partial_{x_{k}}+ A_{k}\right)\right)+q\Big]u(t,x)=0,\ (t,x)\in M_{T}\\
u(0,x)=\phi(x),\ x\in M\\
u(t,x)=f(t,x),\ (t,x)\in \Sigma
\end{cases}
\end{aligned}
\end{align}
where $g^{-1}:=\left(\left( g^{ij}\right)\right)_{1\leq i,j\leq n}$ denote the inverse of metric tensor $g:=\left(\left( g_{ij}\right)\right)_{1\leq i,j\leq n}$,  $\lvert g\rvert=\det(g)$ and the initial value $\phi$ and the Dirichlet data  $f$ are assumed to be non-zero. 
Throughout this article, we denote by $\mathcal{L}_{A,q}$ the following operator
\begin{equation}\label{eq: Operator LAq}
\mathcal{L}_{A,q}:= \partial_{t}-\sum\limits_{j,k=1}^{n}\frac{1}{\sqrt{\lvert g\rvert}}\Big(\partial_{x_{j}}+ A_{j}(t,x)\Big)\left(\sqrt{\lvert g\rvert} g^{jk}\left (\partial_{x_{k}}+ A_{k}(t,x)\right)\right)+q(t,x).
\end{equation}
In this paper, we are interested in determining the convection term $A$ and density coefficient $q$ from the boundary measurements of the solution. To define the boundary operators, we need to have the existence and uniqueness of a solution to the forward problem for IBVP given by  \eqref{eq: Main equation}.

{ Motivated by   \cite{Caro_Kian_Convection_nonlinear,Suman_Manmohan_IPI}, we define the following spaces
\begin{align*}
\begin{aligned}
&	\mathcal{K}_{0}:=\{(f|_{t=0},f|_{\Sigma}):\ f\in H^{1}(0,T,H^{-1}(M))\cap L^{2}(0,T;H^{1}(M))\}\ \ 
\mbox{and} \\
&\mathcal{K}_{T}:=\{(f|_{t=T},f|_{\Sigma}):\ f\in H^{1}(M_T)\}
\end{aligned}
\end{align*}
where we refer to \cite{Evans_Book} for the definition of function spaces $H^{m}(0,T;H^{k}(M)),$  for $k,m\in \mathbb{R}$. Now for $(\phi,f)\in \mathcal{K}_{0}$, it can be shown by following arguments from \cite{Caro_Kian_Convection_nonlinear, Evans_Book, Lion-Magenes2} that there exists a unique solution $\displaystyle u\in H^{1}(0, T, H^{-1}(M))\cap L^{2}(0, T; H^{1}(M))$ of IBVP \eqref{eq: Main equation}.
Based on the existence and uniqueness of solution and following \cite{Caro_Kian_Convection_nonlinear, KU_CMP_Elliptic, Suman_Manmohan_IPI},    we  observe that for any solution $\displaystyle u\in H^{1}(0, T, H^{-1}(M))\cap L^{2}(0, T; H^{1}(M))$ the operator $\displaystyle\mathcal{N}_{A,q}u$ given by  
\begin{align*}
    \begin{aligned}
        \langle \mathcal{N}_{A,q}u,w|_{\partial M^{*}_T}\rangle&:=\int_{M_{T}}\left(-u\partial_t \overline{w}+\langle \nabla_g u,\nabla_g \overline{w}\rangle_g+2u\langle A,\nabla_g \overline{w}\rangle_g +(\delta_g A)u\overline{w}-\lvert A\rvert^2_g u \overline{w}+qu\overline{w}   \right) \ dV_g dt\\
       & \ \ \  -\int_{M}u(0,x)\overline{w}(0,x)dV_{g}
    \end{aligned}
\end{align*}
is well-defined for all $w\in H^1(M_T)$ where $\partial M_{T}^{*}:=(\{T\}\times M)\cup \Sigma.$
Now if we assume the sufficient regularity on the coefficients $A,q$ and Dirichlet data $f$, then as shown in \cite{Caro_Kian_Convection_nonlinear} the operator $\mathcal{N}_{A,q}u$ is given by 
\[\mathcal{N}_{A,q}u:=\Big(u|_{t=T},\Big[\partial_{\nu}u(t,x)+2\nu(x)\cdot A(t,x)u(t,x)\Big]\Big|_{\Sigma}\Big)\]
where $\nu$ stands for the outward unit normal vector to $\partial M$ and $u$ solves the IBVP given by \eqref{eq: Main equation}. This motivates us to define our  input-output operator $\displaystyle \Lambda_{A,q}:\mathcal{K}_{0}\rightarrow \mathcal{K}_{T}^{*}$ by  
\begin{align}\label{eq: Full DN map}
\Lambda_{A,q}(\phi,f):=\mathcal{N}_{A,q}u
\end{align}
where $\mathcal{K}_T^{*}$ stands for dual of $\mathcal{K}_{T}$ and $u$ is solution to the IBVP \eqref{eq: Main equation} when the initial data is $\phi$ and Dirichlet boundary data equal to $f$. }

This work is concerned with the determination of time-dependent coefficients $A$ and $q$ appearing in \eqref{eq: Main equation} using the measurements of the input-output operator $\Lambda_{A,q}$  on a proper subset of $\Sigma$ for the case when $(M,g)$ is an admissible manifold whereby an admissible manifold, we mean the following. 
\begin{definition}(Admissible manifold \cite{Bhattacharyya_IPI,Ferreira_Kenig_Salo_Uhlmann_invention})
We say that a compact Riemannian manifold $(M,g)$ of dimension $n\geq 2$ with boundary $\partial M$, is admissible if $M$ is orientable and $(M,g)$ is a submanifold of $\mathbb{R}\times (\mbox{int}(M_{0}),g_{0})$  where $(M_{0},g_{0})$ is a compact, simply connected Riemannian manifold with boundary $\partial M_{0}$ which is strictly convex in the sense of the second fundamental form and $M_{0}$ has no conjugate points. 
\end{definition}
In order to state the main result of this article, we first need to specify the subset of $\partial M$ where the measurements are given. Now if we write $x\in M$, as  $x:=(x_{1},x')\in \mathbb{R}\times M_{0}$ and  $\varphi(x):=x_{1}$ then  $\partial M$ can be decomposed into the two parts given by 
\[\partial{M}_{+}:=\{x\in \partial M:\ \partial_{\nu}\varphi(x) > 0\}\ \mbox{and}\ \partial{M}_{-}:=\{x\in \partial M:\ \partial_{\nu}\varphi(x)\leq 0\}\] where $\nu(x)$ stands for outward unit normal to $\partial M$ at $x\in \partial M$ and $\partial_{\nu}\varphi$ denote the normal derivative of $\varphi$  with respect to the metric $g$. 
	In this paper, we will be assuming that our boundary measurements are given on slightly bigger than half of $\partial M$. To specify this portion of $\partial M$, we take $\epsilon>0$ small enough and define $\displaystyle \partial M_{\pm,\epsilon/2}$ by 
	\begin{align}\label{Dividing the boundary} 
		\begin{aligned}
		\partial{M}_{+,\epsilon/2}:=\left\{x\in \partial M:\ \partial_{\nu}\varphi(x)\geq \frac{\epsilon}{2}\right\}\ \mbox{and}\ \partial{M}_{-,\epsilon/2}:=\left\{x\in \partial M:\ \partial_{\nu}\varphi(x)< \frac{\epsilon}{2}\right\}\ 
		\end{aligned}
	\end{align} 
 as $\displaystyle \partial M_{-,\epsilon/2}$ is the small enough open neighborhood of $\partial M_{-}$. We denote the  corresponding lateral part of $\Sigma$  by $\Sigma_{+}:= (0,T)\times \partial{M}_{+},$ $\Sigma_{+,\epsilon/2}:=(0,T)\times \partial{M}_{+,\epsilon/2}$, $\Sigma_{-}:=(0,T)\times \partial{M}_{-}$ and $\Sigma_{-,\epsilon/2}:=(0,T)\times \partial{M}_{-,\epsilon/2}$. { We also denote $\displaystyle \partial M^{*}_{T_{\pm}}:= (\{T\}\times M)\cup \Sigma_{\pm} $ and $\displaystyle \partial M^{*}_{T_{\pm,\epsilon/2}}:=  (\{T\}\times M)\cup \Sigma_{\pm,\epsilon/2}.$
 Now, using these notations,
 we define the partial input-output operator by 
	\begin{align}\label{Partial DN map}
		\Lambda^{partial}_{A,q}\big(\phi, f\big)
		:=\mathcal{N}_{A,q}u|_{\partial M^{*}_{T_{-,\epsilon/2}}}. 
	\end{align}	}
Our aim in this article is to recover $A$ and $q$ uniquely from the knowledge of $\Lambda^{partial}_{A,q}$ however, due to gauge invariance, it is impossible to recover these coefficients fully.    Since this is the first work related to the time-dependent convection-diffusion equation on manifolds therefore before stating the main result of this paper,  we first provide quick proof of the gauge invariance associated with our problem. In the Euclidean setting, this has been well observed in prior works; see, for example, \cite{Caro_Kian_Convection_nonlinear, Suman_Manmohan_IPI} and references therein. 
 \begin{definition}(Gauge Invariance)\label{Gauge invariance}
    Let $A^{(i)} \in W^{1,\infty}(M_{T})$ and $q_i \in L^{\infty}(M_T)$ for $i=1,2$. We say  $(A^{(1)}, q_1)$ and $(A^{(2)}, q_2)$ are gauge equivalent if there exists $\Psi \in W_0^{2,\infty}(M_T)$ such that
    \begin{align*}
		A^{(2)}(t,x) = A^{(1)}(t,x) - \nabla_{g}\Psi(t,x)\ \mbox{and} \ q_{2}(t,x) = q_{1}(t,x) - \partial_t\Psi(t,x), \ \mbox{for}\ (t,x)\in M_{T}.
		\end{align*}
 \end{definition}

 \begin{prop}\label{Gauge_equivalent_Prop}
     Suppose $u_1(t,x)$ is a solution to the following IBVP
     \begin{align}\label{eq: u_1 Solution}
     \begin{cases}
     \Big[\partial_{t}-\sum\limits_{j,k=1}^{n}\frac{1}{\sqrt{\lvert g\rvert}}\Big(\partial_{x_{j}}+ A_{j}^{(1)}\Big)\left(\sqrt{\lvert g\rvert} g^{jk}\left (\partial_{x_{k}}+ A_{k}^{(1)}\right)\right)+q_1\Big]u_1(t,x)=0,\ (t,x)\in M_{T}\\
     u_1(0,x)=\phi(x),\ x\in M\\
     u_1(t,x)=f(t,x),\ (t,x)\in \Sigma
\end{cases}
    \end{align}
     and $\Psi \in W_0^{2,\infty}(M_T)$, then $u_2(t,x) = e^{\Psi(t,x)}u_1(t,x)$ satisfies the following IBVP
      \begin{align}\label{u_2 Solution}
     \begin{cases}
\Big[\partial_{t}-\sum\limits_{j,k=1}^{n}\frac{1}{\sqrt{\lvert g\rvert}}\Big(\partial_{x_{j}}+ A_{j}^{(2)}\Big)\left(\sqrt{\lvert g\rvert} g^{jk}\left (\partial_{x_{k}}+ A_{k}^{(2)}\right)\right)+q_2\Big]u_2(t,x)=0,\ (t,x)\in M_{T}\\
u_2(0,x)=\phi(x),\ x\in M\\
u_2(t,x)=f(t,x),\ (t,x)\in \Sigma
\end{cases}
     \end{align}
   where $A^{(2)}(t,x) = A^{(1)}(t,x) - \nabla_{g}\Psi(t,x)\ \mbox{and} \ q_{2}(t,x) = q_{1}(t,x) - \partial_t\Psi(t,x).$  Now if  $\Lambda_{A^{(i)}, q_i}$ for $i=1,2$, are the input-output operators associated with $u_i$ and defined by \eqref{eq: Full DN map} then 
     \begin{align*}
         \Lambda_{A^{(1)}, q_1}(\phi,f) = \Lambda_{A^{(2)}, q_2}(\phi,f),\ \mbox{for all}\ (\phi,f)\in \mathcal{K}_{0}.
     \end{align*} 
 \end{prop}
 \begin{proof}
    Substituting $u_1(t,x) = e^{-\Psi(t,x)}u_2(t,x)$ in Equation \eqref{eq: u_1 Solution}, from simple computations, we get
    \begin{align*}
       0= \mathcal{L}_{A^{(1)},q_1}(e^{-\Psi(t,x)}u_2(t,x)) =  \mathcal{L}_{A^{(1)} - \nabla_{g}\Psi,q_1- \partial_t\Psi}u_2(t,x) = \mathcal{L}_{A^{(2)},q_2}u_2(t,x),\ \quad (t,x) \in M_T
    \end{align*}
    and
    \begin{align*}
       & u_2(0,x) = e^{\Psi(0,x)}u_1(0,x) = \phi(x), \quad x \in M,\\
      & u_2(t,x) = e^{\Psi(t,x)}u_1(t,x) = f(t,x), \quad (t,x) \in \Sigma.
    \end{align*}
    Hence $u_2$ solves \eqref{u_2 Solution}. Also, we have that 
    \begin{align*}
    \begin{aligned}
       & u_2(T,x) = e^{\Psi(T,x)}u_1(T,x) = u_1(T,x), \quad x \in M,\ \  \partial_{\nu}u_2\Big|_{\Sigma} = \left(e^{\Psi}(\partial_{\nu}\Psi u_1 + \partial_{\nu}u_1\right)\Big|_{\Sigma} = \partial_{\nu}u_1\Big|_{\Sigma},\\
    \mbox{and}  & \left(\nu\cdot A^{(2)}u_2\right)\Big|_{\Sigma} = \left(\nu\cdot (A^{(1)} - \nabla_g\Psi(e^{\Psi}u_1)\right)\Big|_{\Sigma}  = \left(\nu \cdot A^{(1)}u_1\right)\Big|_{\Sigma} 
    \end{aligned}
    \end{align*}
    where in the above equations, we have used the fact that $\Psi \in W_0^{2,\infty}(M_T)$. Thus combining the above equations together with  \eqref{eq: Full DN map} we get \begin{align*}
         \Lambda_{A^{(1)}, q_1}(\phi,f) = \Lambda_{A^{(2)}, q_2}(\phi,f),\ \mbox{for all}\ (\phi,f)\in \mathcal{K}_{0}.
     \end{align*} 
 \end{proof}
	\noindent {With this preparation, we are ready to }state the main result of this paper as follows. 
	\begin{theorem}\label{th:main theorem}
		Let $(M,g)$ be an admissible manifold. Let $\displaystyle A^{(i)}\in W^{1,\infty}(M_{T};T^{*}(M))$ for $i=1,2$ given by $\displaystyle A^{(i)}(t,x)=\sum\limits_{j=1}^{n}A^{(i)}_{j}(t,x)dx^{j},$ in local coordinates on $(M,g)$ and $q_{i}\in L^{\infty}(M_T)$ for $i=1,2$. Suppose $u_{i}$ for $i=1,2,$ is solution to \eqref{eq: Main equation} when $(A,q)=(A^{(i)},q_{i})$ for $i=1,2$ and $\Lambda^{partial}_{A^{(i)},q_{i}}$ are input-output operator given by \eqref{Partial DN map} corresponding to $u_{i}$ for $i=1,2$. Now for $\epsilon>0$ small enough if 
		\begin{align}
  \label{eq: Equal DN map }
	\Lambda^{partial}_{A^{(1)},q_{1}}\big(\phi,f\big) = \Lambda^{partial}_{A^{(2)},q_{2}}\big(\phi,f\big),\ \mbox{for all}\ (\phi,f)\in \mathcal{K}_{0}
		\end{align}  
		then there exists a function $\Psi\in W_0^{2,\infty}(M_{T})$ such that 
		\begin{align*}
			A^{(1)}(t,x)-A^{(2)}(t,x)=\nabla_{g}\Psi(t,x)\ \mbox{and} \ q_{1}(t,x) - q_{2}(t,x) = \partial_t\Psi(t,x), \ \mbox{for}\ (t,x)\in M_{T}
		\end{align*}
		provided $A^{(1)}(t,x)=A^{(2)}(t,x)$, for $(t,x)\in \Sigma.$
	\end{theorem}
\begin{remark}
\begin{enumerate}
\item { Observe that the measurement data used in Theorem \ref{th:main theorem} is an input-output map, unlike the usual Dirichlet to Neumann (DN) map used in the Euclidean setting. This is due to the fact that in our boundary Carleman estimate stated in Theorem \ref{Boundary Carleman estimate Theorem}, we do not have the estimate on weighted  $L^2$ norm of solution $u$ at $t=T$ which is because of the choice of  weight function $\varphi(t,x):=\lambda^2\beta^2 t+\lambda x_1$ where $\beta\in (0,1)$ while in Euclidean setting one can actually take $\beta=1$ which helped one to the  Carleman estimate with a bound on the weighted  $L^2$ norm of the solution $u$ at $t=T$. We refer to \cite{Caro_Kian_Convection_nonlinear, Soumen_Manmohan_EECT} for details about it. However, if we assume the coefficients are small enough then we can obtain the boundary Carleman estimate (see Theorem $3.1$ in \cite{Suman_Manmohan_IPI}) with a bound on the weighted  $L^2$ norm of the solution $u$ at $t=T$  and can determine the coefficients from the knowledge of DN map measured on a suitable subset of $\Sigma.$}
\item We observe that because of gauge invariance proved in Proposition \ref{Gauge_equivalent_Prop}, it is impossible to prove that  $A^{(1)}=A^{(2)}$ and $q_1=q_2$ in $M_{T}$ from the given hypothesis of Theorem \ref{th:main theorem}. As far as the uniqueness issue is concerned,  gauge invariance {guarantees} that {the result} obtained in Theorem  \ref{th:main theorem} is optimal.
\item Unique recovery of $A$ and $q$ is also possible with some extra conditions on the unknown vector field $A$. For instance, if $A$ is divergence-free, that is,  $ \delta_g A =0$ in $M_{T}$, then one can recover both $A$ and $q$ uniquely in $M_{T}$. This divergence-free condition has been exploited in earlier works to get the full recovery,  please refer \cite{Suman_Manmohan_IPI, Soumen_Manmohan_EECT}. 
\item On the other side, if we assume that the vector field $A$ is time-independent, then recovery of $A$ is possible up to a potential of the form $\nabla_g\Psi(x)$ and $q$ can be fully recovered in $M_{T}$. The proof follows similarly as we have done for Theorem \ref{th:main theorem}. 
\end{enumerate}
\end{remark}
The problem considered in this article can be put under the umbrella of Calder\'on type inverse problems for parabolic partial differential equations (PDEs), which was initially proposed by Calder\'on in \cite{Calderon_problem} for elliptic PDEs and studied by Nachman \cite{Nachman_Calderon-problem} in two dimensions and by Sylvester-Uhlmann in \cite{Sylvester_Uhlmann_Calderon_problem_1987} in dimension three and higher. Analogous {problems} for parabolic and hyperbolic PDEs have been studied in \cite{Avdonin_Seidman_parabolic, Isakov_Completeness, Nakamura_Sasayama_Parabolic, Rakesh_Symes}. Choulli-Kian in \cite{Choulli_Kian_stability_parabolic_MCRF} derived a stability estimate for recovering the time-dependent coefficient, which is a product of functions depending only on time and only space variables, from the boundary measurements.  We also refer to \cite{Choulli_Abstract_IP_Applications}, where an abstract inverse problem for parabolic pde is studied. All these works are concerned with the recovery of zeroth order perturbation of elliptic, parabolic, and hyperbolic PDEs from full boundary data. In \cite{Choullin_Kian_JMPA},  the recovery of general time-dependent zeroth order perturbation of heat operator from partial boundary measurements is considered. Inverse problems of recovering the coefficients appearing in the steady state convection-diffusion from full and partial boundary measurements in Euclidean geometry have been studied in  \cite{Brown_Salo_static_convection-diffusion, Cheng_Yamamoto_Global_Convection_2D_DN, Cheng_Yamamoto_Steady, Christofol_Soccorsi_magnetic, Ferriera_Kenig_Sjostrand_Uhlmann_magnetic, Kim_Salo_Static convection IPI, Nakamura_Sun_Uhlmann_Magnetic, Pohjola_steady_state_CD, Sun_magnetic}. Recovery of first-order perturbation of a parabolic pde from final and single measurement has been studied in \cite{Deng_Yu_Yang_First_order_parabolic_1995} and \cite{Chen_Yamamoto_parabolic} respectively.  In \cite{Bellassoued_Rassas_Convection-diffusion_JIIP}, stable recovery of time-dependent coefficients appearing in a convection-diffusion from full boundary data has been studied. Choulli-Kian in \cite{Choullin_Kian_JMPA} proved a stability estimate for recovering a time-dependent potential from partial boundary data, and motivated by their work, authors of \cite{Suman_Manmohan_IPI} proved the unique recovery of time-dependent {coefficients} appearing in a convection-diffusion equation from partial boundary data. In \cite{Suman_Manmohan_IPI}, a uniqueness result is proved with a smallness assumption on the convection term, which {is later} on removed in a recent work of  \cite{Soumen_Manmohan_EECT} where stability estimates for recovering the time-dependent coefficients of convection-diffusion equation from partial boundary data are derived. We also refer to \cite{Bellassoued_Rassas_Convection-diffusion_JIIP} and  \cite{Bellassoued_Oumaima_Convection-diffusion_IP}  where stability estimates for convection-diffusion equation from full and partial boundary data is studied, respectively. Recently, in \cite{Caro_Kian_Convection_nonlinear, Feizmohammadi_Kian_Uhlmaann_quasilinear convection-diffusion}, inverse problems related to nonlinear convection-diffusion equation is studied.
    In all the above-mentioned works, the inverse problems of recovering coefficients appearing in parabolic PDEs from full and partial boundary measurements in Euclidean geometry are considered. The inverse problems related to steady-state convection-diffusion equation in Riemannian geometry are considered in prior works (see, for example \cite{Bhattacharyya_IPI, Ferreira_Kenig_Salo_Uhlmann_invention, Kenig_Salo_Analysis and PDEs, Krupchyk_Uhlmann_CPDE_2017}) however the recovery of time-dependent coefficients appearing in parabolic PDEs in Riemannian geometry has {not been} considered in prior works, and this is the main objective of this paper. To the best of our knowledge, this is the first work {that considers} the partial data inverse problem for recovering both first and zeroth order time-dependent perturbations of evolution equations in Riemannian geometry. Next, we mention works on inverse problems related to hyperbolic and dynamical Schr\"odinger equations, which are closely related to the study of this work. 
	Inspired by \cite{Bukhgeim_Uhlmann_paartial data}  and \cite{Bellassoued_Jellali_Yamamoto_Lipschitz_stability_hyperbolic, Bellassoued_Jellali_Yamamoto_stability_hyperbolic} authors of  \cite{Kian_Damping_partial_data_2016, Kian_Uniqueness_partial_data_2016, Manmohan_Venky, Mishra_Vashisth_Uniqueness_AA} studied the unique recovery of time-dependent coefficients appearing in a hyperbolic pde from partial boundary measurements. 
	We also refer to \cite{Ibtissem_dynamical Schroedinger, Bellassoued_Ferreira_anisotropic_Schrodinger, Bellassoued_Ferreira_anisotropic_wave, Bellassoued_Rezig_dynamical Schroedinger_MCRF, Bellassoued_Rezig_Hyperbolic_Riemannian_AGAG, AlI_Ilmavirta_Kian_Oksanen_wave-manifold, Kian_Soccorsi_Holder_stability_Scrodinger, Kian_Tetlow_Holder_stability_dynamical Scrodinger, Salazar_time-dependent_first_order_perturbation_2013} for inverse problems related to hyperbolic PDEs and dynamical Schr\"odinger equation in Euclidean as well as in Riemannian geometry, where time-dependent coefficients are recovered from full boundary measurements. 
	
	The rest of the article is organized as follows. In section  \ref{Boundary Carleman estimate section}, we derive the boundary and interior Carleman estimates {which we will use in section \ref{Construction of GO solutions} to construct} the exponentially growing as well as decaying solutions. {The main Theorem \ref{th:main theorem} of the article  will be proved in section \ref{Proof of main theorem}}. Finally, we conclude the article with Acknowledgements.
	\section{Boundary and interior Carleman Estimates}\label{Boundary Carleman estimate section}
	The present section is devoted to prove a boundary and interior Carleman estimates. The boundary Carleman estimate will be required to estimate the boundary terms in the integral identity where there is no measurement, and the interior Carleman estimate will be required to construct the geometric optics solutions (GO) for $\mathcal{L}_{A,q}$ and its formal $L^{2}-$-adjoint. 
	For a compact Riemannian manifold $(M,g)$ with boundary denoted by $\partial M$, we denote by $dV_{g}$ the volume form on $(M,g)$ and by  $dS_{g}$ the induced volume form on $\partial M$. Then the $L^{2}$-norm of a function $u$ on $M$ and $f$ on $\partial M$ are given by 
	\[\lVert u\rVert_{L^{2}(M)}:=\left(\int_{M}\lvert u(x)\rvert^{2}\ dV_{g}\right)^{1/2}\ \mbox{and}\ \lVert f\rVert_{L^{2}(\partial M)}:=\left(\int_{\partial M}\lvert f(x)\rvert^{2}\ dS_{g}\right)^{1/2}, \ \mbox{respectively}.\]
	We denote by $L^{2}(M)$ as the space of all functions $u$ defined on $M$  for which $\lVert u\rVert_{L^{2}(M)}<\infty$ and $L^{2}(\partial M)$ as the space of all functions $f$ defined on $\partial M$ for which $\lVert f\rVert_{L^{2}(\partial M)}<\infty.$ Then  $\left(L^{2}(M),\lVert\cdot\rVert_{L^{2}(M)}\right)$ and $\left(L^{2}(\partial M),\lVert\cdot\rVert_{L^{2}(\partial M)}\right)$ are Hilbert spaces with respect to the inner-products defined by 
	\begin{align*}
		\langle f,g\rangle_{L^{2}(M)}:=\int_{M}f(x)\overline{g(x)}\ dV_{g}
	\end{align*}
	and \[ \langle f,g\rangle_{L^{2}(\partial M)}:=\int_{\partial M}f(x)\overline{g(x)}\ dS_{g}\] respectively. \\
	
	We state the boundary Carleman estimate as follows. 
	\begin{theorem}\label{Boundary Carleman estimate Theorem}
		Let $(M,g)$ be an admissible manifold. 	For $\displaystyle \beta\in \left(\frac{1}{\sqrt{3}},1\right)$, let  $\varphi(t,x):=\lambda^{2}\beta^{2}t+\lambda x_{1}$, $A \in \left(W^{1,\infty}(M_{T})\right)^{n}$ and $q\in L^{\infty}(M_T)$. Then there exists a constant $C>0$ depending only on $M, T$, $A$ and $q$ such that 
		\begin{align}\label{eq: BC estimate}
			\begin{aligned}
				&\lambda^2  \|e^{-\varphi} u\|_{L^2(M_T)}^2 +\left\|e^{-\varphi} \nabla_g u\right\|_{L^2(M_T)}^2+ \|e^{-\varphi(T, \cdot)} \nabla_g u(T, \cdot)\|_{L^2(M)}^2   +\lambda\int_{ \Sigma_{+}} e^{-2\varphi}\lvert  \partial_{\nu}u(t,x)\rvert^{2}\langle \nu,e_{1}\rangle_{g}\ dS_{g}dt\\
				&\ \ \leq C\lb \| e^{-\varphi} \mathcal{L}_{A,q} u\|_{L^2(M_T)}^2 +\lambda^2  \|e^{-\varphi(T, \cdot)}u(T, \cdot )\|_{L^2(M)}^2 +\lambda \int_{ \Sigma_{-}} e^{-2\varphi}\lvert  \partial_{\nu}u(t,x)\rvert^{2}|\langle \nu,e_{1}\rangle_{g}| \ dS_{g}dt\rb
			\end{aligned}
		\end{align}
		hold for $\lambda$ large enough and for all $u\in C^{2}(\overline{M_{T}})$ satisfying the following
		\begin{align*}
			\begin{aligned}
				u(0,x)=0, \ \mbox{for} \ x\in M\ \ \mbox{and}\ u(t,x)=0,\ \mbox{for} \ (t,x)\in \Sigma. 
			\end{aligned}
		\end{align*}

		\begin{proof}
			In order to prove the weighted  $H^1-L^2$  estimate given by \eqref{eq: BC estimate}, we need to convexify the Carleman weight $\varphi$. This convexification will help us to absorb the first-order perturbation  $A$ appearing in  $\mathcal{L}_{A,q}$, which has been used in \cite{Caro_Kian_Convection_nonlinear, Soumen_Manmohan_EECT} for Euclidean case and \cite{Ferreira_Kenig_Salo_Uhlmann_invention} for anisotropic magnetic Schr\"odinger operator. Now  for $s>0$, we denote the convexified weight function by $\varphi_{s}$ and define by  
			\begin{align}\label{eq: Convexified Carleman weight} 
				\varphi_s(t, x) := \varphi(t, x) - \frac{s(x_1 + 2 \ell)^2}{2}=\lambda^{2}\beta^{2}t+\lambda x_{1}- \frac{s(x_1 + 2 \ell)^2}{2}.
			\end{align}
			A direct computation gives 
			\begin{align}\label{Calculation with phis}
				\partial_t\vph_s=  \lambda^2 \beta^2, \ \  \partial_{x_1} \vph_s =  \lambda  - s (x_1 +2 \ell),\ \ 
				\partial_{x_1}^2 \vph_s =  -s ,\ \mbox{and}\  \lvert\partial_{x_1} \vph_s \rvert^2 =  \lb\lambda  - s (x_1 +2 \ell)\rb^2.
			\end{align}
			Before we proceed further, let us observe that 
			\begin{align*}
				\mathcal{L}_{A,q}v(t,x)=\partial_{t}v(t,x)-\partial_{x_{1}}^{2}v(t,x)-\Delta_{g_0}v(t,x)-2\langle A(t,x),\nabla_{g}v(t,x)\rangle_{g}+\tilde{q}(t,x) v(t,x)
			\end{align*}
			where $\langle\cdot,\cdot\rangle_{g}$ and $\nabla_{g}$ denote the inner-product and gradient operator  w.r.t. metric $g$ respectively and $\tilde{q}(t, x) := q(t, x) -\delta_g A(t, x) - \lvert A(t,x)\rvert_g^2$, here in expression of $\tilde{q}$, $\delta_{g}A$ given by  $$\delta_g A  = \frac{1}{\sqrt{\lvert g\rvert}}\sum_{j, k = 1}^n \partial_j \left(g^{jk}\sqrt{|g|} A_k\right)$$ is known as the divergence operator w.r.t. to metric $g$  and $\lvert A\rvert_{g}^{2}=\sum_{j,k=1}^{n}g^{jk}A_{j}A_{k}$. With this, we define the conjugated operator $P_{s}$ with a convexified weight function $\varphi_{s}$ by
			\begin{align}\label{Conjugated operator}
				P_{s}v:=e^{-\varphi_s}\Lc_{A, q}\lb e^{\varphi_s}v\rb=e^{-\varphi_s}\lb \partial_{t}-\partial_{x_{1}}^{2}-\Delta_{g_0}-2\langle A,\nabla_{g}\rangle_{g}+\tilde{q}\rb \lb e^{\varphi_s}v\rb. 
			\end{align} 
			Upon expanding the above expression,  $P_{s}$ will take the following form 
			\begin{align*}
				P_{s}v(t,x)
				&\ =\lb \partial_t  + (\partial_t\vph_s)\rb v(t,x)-\lb\partial^2_{x_1} + 2\partial_{x_1}\vph_s\partial_{x_1}+(\partial_{x_1}\vph_s)^2+ \partial_{x_1}^2\vph_s \rb v -\Delta_{g_0} v(t,x)\\
				&\qquad - 2 \langle A (t,x), \nabla_g v(t,x)\rangle_g - 2 \langle A(t,x) , \nabla_g \vph_{s}(t,x) \rangle_g v(t,x)  + \tilde{q}(t,x)v(t,x).
			\end{align*}
			Using \eqref{Calculation with phis} and the fact that $(M,g)$ is admissible manifold, we get 
			\begin{align*}
				\begin{aligned}
					P_{s}v(t,x)&=\partial_{t}v(t,x)+\lambda^{2}\beta^{2}v(t,x)-\left( \partial_{x_{1}}^{2}+2(\lambda-s(x_{1}+2\ell))\partial_{x_{1}}+(\lambda-s(x_{1}+2\ell))^{2}-s\right)v(t,x)\\
					&\qquad -\Delta_{g_0}v(t,x)	- 2 \langle A (t,x), \nabla_g v(t,x)\rangle_g-2  \lb\lambda - s(x_1 + 2\ell)\rb g^{1k}A_k(t,x) v(t,x)+\tilde{q}(t,x)v(t,x). 
				\end{aligned}
			\end{align*}
			Now if we define $P_1,P_2$ and $P_3$ by
			\begin{align*}
				\begin{aligned}
					P_1 v(t,x)&:= \left(\partial_{t}v - 2(\lambda- s (x_1 + 2 \ell))\partial_{x_1}v + 4sv\right)(t,x),\\
					P_2 v(t,x)&:= \left (-\partial^2_{x_1}v - \Delta_{g_0} v - \lambda^2(1 - \beta^2)v + 2 \lambda s(x_1 + 2 \ell) v - s^2 (x_1 + 2 \ell)^2 v -  3sv\right)(t,x)\\
					&:=\lb -\partial^2_{x_1} - \Delta_{g_0}+\mathcal{K}(x_{1})\rb v(t,x), \ \mbox{where} \ \mathcal{K}(x_{1}):=2 \lambda s(x_1 + 2 \ell)- \lambda^2(1 - \beta^2) -  s^2 (x_1 + 2 \ell)^2-3s\\
					P_3 v(t,x)&:= -2 \langle A(t,x) , \nabla_g v(t,x)\rangle_{g} -2 \lb \lambda - s(x_1 + 2\ell) \rb g^{1k}A_k(t,x) v(t,x) + \tilde{q}(t,x)v(t,x)
				\end{aligned}
			\end{align*}
			then one can check that $P_{s}v(t,x)$ has the following compact form 
			\begin{align}\label{Expression for Psv}
				P_{s}v(t,x)=P_{1}v(t,x)+P_{2}v(t,x)+P_{3}v(t,x).
			\end{align}
			Our first aim is to estimate the $L^{2}$ norm of $P_{s}v$ on $M_{T}$, therefore we define  $I_{s}$ by 
			\begin{align*}
				\begin{aligned}
					I_s &:=\int_{M_T}\lvert P_{s}v(t,x)\rvert^{2} \ dV_{g}(x)dt
					=\int_{M_T}\lvert P_{1}v(t,x)+P_{2}v(t,x) + P_{3}v(t,x)\rvert^{2} \ dV_{g}dt\\
					&\geq \frac{1}{2} \int_{M_T}\left( P_{1}v(t,x)+ P_{2}v(t,x)\right)^{2} \ dV_{g}dt - \int_{M_T} \lvert P_{3}v(t,x)\rvert^2\ dV_{g}dt\\
					&\geq  \int_{M_T}P_{1}v(t,x) P_{2}v(t,x)\ dV_{g}dt - \int_{M_T} |P_{3}v(t,x)|^2\ dV_{g}dt.
				\end{aligned}
			\end{align*}
			This gives us 
			\begin{align}\label{First estimate on Is}
				\begin{aligned}
					I_{s}\geq \underbrace{\int_{M_T}P_{1}v(t,x) P_{2}v(t,x)\ dV_{g}dt}_{I_{s, 1}} - \underbrace{\int_{M_T} |P_{3}v(t,x)|^2\ dV_{g}dt}_{I_{s, 2}}.
				\end{aligned}
			\end{align}
			We aim to estimate the right-hand side of  \eqref{First estimate on Is}. To do that, we start with the first term in the above inequality and, therefore consider
			\begin{align*}
				\begin{aligned}
					P_{1}v(t, x)P_{2}v(t,x)&=- \partial_{t}v(t,x)\left(\partial^2_{x_1}v+\Delta_{g_{0}} v\right )(t,x) +  \mathcal{K}(x_{1})v(t,x)\partial_{t}v(t,x)-4sv(t,x)\left (\partial_{x_{1}}^{2}+\Delta_{g_{0}}\right)v(t,x) \\
					&\quad +4s\mathcal{K}(x_{1})\lvert v(t,x)\rvert^{2}
					+ 2(\lambda- s (x_1 + 2 \ell))\partial_{x_1}v(t,x)\left (\partial^2_{x_1}v+ \Delta_{g_{0}} v\right )(t,x)\\
					&\qquad  -2\mathcal{K}(x_{1})(\lambda- s (x_1 + 2 \ell))v(t,x)\partial_{x_1}v(t,x). 
				\end{aligned}
			\end{align*}
			Now consider $I_{s,1}$  from \eqref{First estimate on Is}
			\begin{align*}
				\begin{aligned}
					I_{s, 1} &= \int_{M_T}P_{1}v(t,x)P_{2}v(t,x)\ dV_{g}dt =-\int_{M_T} \partial_{t}v(t,x)\left(\partial^2_{x_1}v+\Delta_{g_{0}} v\right )(t,x)\ dV_{g}dt \\ 
					&\quad + 4s\int_{M_T}\mathcal{K}(x_{1})\lvert v(t,x)\rvert^{2}\ dV_{g}dt + \frac{1}{2}\int_{M_T}\mathcal{K}(x_{1})\partial_{t}\lvert v(t,x)\rvert^{2}\ dV_{g}dt\\
					&\quad - 4s\int_{M_T}v(t,x)\left (\partial_{x_{1}}^{2}v+\Delta_{g_{0}}v\right)(t,x)\ dV_{g}dt\\
					&\ \ \ \  + 2\int_{M_T}(\lambda-s (x_1 + 2 \ell))\partial_{x_1}v(t,x)\left (\partial^2_{x_1}v+ \Delta_{g_{0}} v\right)(t,x)\ dV_{g}dt\\
					&\ \ \ \ -\int_{M_T}\mathcal{K}(x_{1})(\lambda- s (x_1 + 2 \ell))\partial_{x_1}\lvert v(t,x)\rvert^{2}\ dV_{g}dt\\
					& := I_{1}+I_{2}+I_{3}+I_{4}+I_{5}+I_{6}
				\end{aligned}
			\end{align*}
			where 
			\begin{align*}
				\begin{aligned}
					&I_{1}:=-\int_{M_T}\partial_{t}v(t,x)\left(\partial^2_{x_1}v+\Delta_{g_{0}} v\right )(t,x)\ dV_{g}dt;\ \ I_{2}:= 4s\int_{M_T}\mathcal{K}(x_{1})\lvert v(t,x)\rvert^{2} \ dV_{g}dt\\
					&I_{3}:=\frac{1}{2}\int_{M_T}\mathcal{K}(x_{1})\partial_{t}\lvert v(t,x)\rvert^{2}\ dV_{g}dt;\ \ I_{4}:= -4s\int_{M_T}v(t,x)\left (\partial_{x_{1}}^{2}v+\Delta_{g_{0}}v\right)(t,x)\ dV_{g}dt\\
					&I_{5}:= 2\int_{M_T}(\lambda-s (x_1 + 2 \ell))\partial_{x_1}v(t,x)\left (\partial^2_{x_1}v+ \Delta_{g_{0}} v\right)(t,x)\ dV_{g}dt\\
					&I_{6}:=-\int_{M_T}\mathcal{K}(x_{1})(\lambda- s (x_1 + 2 \ell))\partial_{x_1}\lvert v(t,x)\rvert^{2}\ dV_{g}dt.
				\end{aligned}
			\end{align*}
			In order to estimate $I_{s,1}$ in \eqref{First estimate on Is}, we need to estimate each  $I_j$ for $1\leq j\leq 6$. To estimate these $I_j's$, we use integration by parts repeatedly along with initial and boundary conditions on $v$. Consider 
			\begin{align}
				\label{eq: I1}  I_{1} &=-\int_{M_T} \partial_{t}v(t,x)\left(\partial^2_{x_1}v+\Delta_{g_{0}} v\right )(t,x)\ dV_{g}dt = \frac{1}{2} \int_M \left| \nabla_g v (T, x)\right|_g^2 dV_g.
			\end{align}
			\begin{align}\label{eq : I2}
				\begin{aligned}
					I_{2} &= 4s\int_{M_T}\mathcal{K}(x_{1})\lvert v(t,x)\rvert^{2}\ dV_{g}dt\\
					&=-4s\int_{M_T}\left (\lambda^2(1 - \beta^2) - 2 \lambda s(x_1 + 2 \ell)  + s^2 (x_1 + 2 \ell)^2+3s\right)\lvert v(t,x)\rvert^{2}\ dV_{g}dt.
				\end{aligned}
			\end{align}
			\begin{align*}
				\begin{aligned} I_{3}&=\frac{1}{2}\int_{M_T}\mathcal{K}(x_{1})\partial_{t}\lvert v(t,x)\rvert^{2}\ dV_{g}dt = \frac{1}{2} \int_{M} \mathcal{K}(x_1) |v(T, x)|^2 dV_{g}\\
					&=\frac{1}{2}  \int_{M}\left(-\lambda^{2}(1-\beta^{2})+2\lambda s(x_{1}+2\ell)-s^{2}(x_{1}+2\ell)^{2}-3s\right) \lvert v(T,x)\rvert^{2}\ dV_{g}.
				\end{aligned}
			\end{align*}
			Recall $\ell \leq (x_1 +2\ell) \leq 3\ell$, therefore choosing $\lambda$ large enough, we obtain
			\begin{align}\label{eq: I3}
				\begin{aligned}
					I_3  &\geq -C\lambda^2 \|v(T, \cdot)\|_{L^2(M)}^2, \ \mbox{ for some constant  } C>0\ \mbox{independent of $\lambda$}.
				\end{aligned}
			\end{align}
			Next, consider
			\begin{align}\label{eq: I4}  
				I_{4} &= -4s\int_{M_T}v(t,x)\left (\partial_{x_{1}}^{2}v+\Delta_{g_{0}}v\right)(t,x)\ dV_{g}dt= 4s \int_{M_T} \left| \nabla_g v (t, x)\right|_g^2 dV_gdt.
			\end{align}
			The next integral in the line is  
			\begin{align*}
				\begin{aligned}
					I_{5}&=2\int_{M_T}(\lambda-s (x_1 + 2 \ell))\partial_{x_1}v(t,x)\left (\partial^2_{x_1}v+ \Delta_{g_{0}} v\right)(t,x)\ dV_{g}dt\\
					&=2\int_{M_T}(\lambda-s (x_1 + 2 \ell))\partial_{x_{1}}v(t,x) \Delta_{g}v(t,x)\ dV_{g}dt.
				\end{aligned}
			\end{align*}
			Using the integration by parts, we have 
			\begin{align*}
				\begin{aligned}
					I_{5}&=-2\int_{M_T}\partial_{x_{1}}v(t,x)\Big\langle\nabla_{g}v,\nabla_{g}\left (\lambda -s(x_{1}+2\ell)\right)\Big\rangle_{g}\ dV_{g}dt \\
					&\qquad - 2\int_{M_T}(\lambda-s (x_1 + 2 \ell))\Big\langle \nabla_{g}v,\partial_{x_{1}}\nabla_{g}v\Big\rangle_{g}\ dV_{g}dt\\
					& \qquad \quad +2\int_{\Sigma}(\lambda-s (x_1 + 2 \ell))\partial_{x_{1}}v(t,x) \partial_{\nu}v(t,x)\ dS_{g}dt
				\end{aligned}
			\end{align*} 
			where $\nu(x)$ is outward unit normal vector to $\partial M$ at $x\in \partial M$, $\partial_{\nu}v(t,x)$ stands for the normal derivative w.r.t. $x$ of $v$ at $(t,x)\in (0,T)\times\partial M$ and $dS_{g}$ denote the surface measure on $\partial M$. Again using the integration by parts, we have that 
			\begin{align*}
				\begin{aligned}
					I_{5}&=2s\int_{M_T}\lvert \partial_{x_{1}}v(t,x)\rvert^{2}\ dV_{g}dt  -\int_{\Sigma} (\lambda-s (x_1 + 2 \ell))\langle \nu,e_{1}\rangle_{g}\lvert \nabla_{g}v(t,x)\rvert_{g}^{2}\ dS_{g}dt\\
					&\qquad -s\int_{M_T} \lvert \nabla_{g}v(t,x)\rvert_{g}^{2}\ dV_{g}dt  +2\int_{\Sigma}(\lambda-s (x_1 + 2 \ell))\partial_{x_{1}}v(t,x) \partial_{\nu}v(t,x)\ dS_{g}dt. \end{aligned}
			\end{align*}
			Now $v|_{(0,T)\times \partial M}=0,$ implies that  $\displaystyle\nabla_{g}v|_{\Sigma}=\left (\partial_{\nu}v\right)\nu$ and $\displaystyle\partial_{x_{1}}v|_{\Sigma}=\langle \nabla_{g}v,e_{1}\rangle_{g}|_{\Sigma}=\partial_{\nu}v\langle \nu,e_1\rangle_{g}$. Using these, we get 
			\begin{align*}
				I_{5} =s\int_{M_T}\left(\lvert \partial_{x_{1}}v(t,x)\rvert^{2}-\lvert \nabla_{g_{0}}v(t,x)\rvert_{g_{0}}^{2}\right) dV_{g}dt +\int_{\Sigma}(\lambda-s (x_1 + 2 \ell))\lvert \partial_{\nu}v(t,x)\rvert_{g}^{2}\langle \nu,e_{1}\rangle_{g}\ dS_{g}dt.
			\end{align*}
			Combining $I_5$  and $I_4$, we have
			\begin{align}\label{eq: I4 and I5}
				\begin{aligned}
					I_4 + I_{5}& =5 s\int_{M_T}\lvert \partial_{x_{1}}v(t,x)\rvert^{2} dV_{g}dt+3s \int_{M_T} \lvert \nabla_{g_{0}}v(t,x)\rvert_{g_{0}}^{2} dV_{g}dt \\
					&\qquad +\int_{\Sigma}(\lambda-s (x_1 + 2 \ell))\lvert \partial_{\nu}v(t,x)\rvert_{g}^{2}\langle \nu,e_{1}\rangle_{g}\ dS_{g}dt\\
					&\geq 3 s  \left\| \nabla_g v\right\|_{L^2(M_T)}^2 + \int_{\Sigma} (\lambda-s (x_1 + 2 \ell))\lvert \partial_{\nu}v(t,x)\rvert_{g}^{2}\langle \nu,e_{1}\rangle_{g}\ dS_{g}dt.
				\end{aligned}
			\end{align}
			Next, we consider the  last term of $I_{s, 1}$
			\begin{align*}
				\begin{aligned}
					I_{6}&=-\int_{M_T}\mathcal{K}(x_{1})(\lambda- s (x_1 + 2 \ell))\partial_{x_1}\lvert v(t,x)\rvert^{2}\ dV_{g}dt\\
					&=\int_{M_T}\left(\lambda -s(x_{1}+2\ell)\right)\lvert v(t,x)\rvert^{2}\partial_{x_{1}}\mathcal{K}(x_{1})\ dV_{g}dt-s\int_{M_T}\mathcal{K}(x_{1}) \lvert v(t,x)\rvert^{2}\ dV_{g}dt\\
					&=2s\int_{M_T}\left ( \lambda -s(x_{1}+2\ell)\right)^{2} \lvert v(t,x)\rvert^{2}\ dV_{g}dt\\
					&\ \ \ +s\int_{M_T}\left (\lambda^2(1 - \beta^2) - 2 \lambda s(x_1 + 2 \ell)  + s^2 (x_1 + 2 \ell)^2+3s \right)\lvert v(t,x)\rvert^{2}\ dV_{g}dt. 
				\end{aligned}
			\end{align*}
			After simplifying, we get 
			\begin{align}
				\begin{aligned}
					I_{6}&=s\lambda^{2}(3-\beta^{2})\int_{M_T}\lvert v(t,x)\rvert^{2}_{g}\ dV_{g}dt-6\lambda s^{2}\int_{M_T}(x_{1}+2\ell)\lvert v(t,x)\rvert_{g}^{2}\ dV_{g}dt\\
					&\ \ \ \ \ +3s^{2}\int_{M_T}\lb s(x_{1}+2\ell)^{2} +1\rb \lvert v(t,x)\rvert_{g}^{2}\ dV_{g}dt.
				\end{aligned}
			\end{align}
			Next, we estimate $I_{s, 2}$ in the following way:
			\begin{align}\label{eq: IS2}
				I_{s, 2} &= \int_{M_T} \lvert P_{3}v(t,x)\rvert^2 dV_{g}dt\nonumber\\
				&= \int_{M_T} \left|-2 \langle A(t,x) , \nabla_g v(t,x)\rangle_{g} - 2 \lb \lambda - s(x_1 + 2\ell)\rb g^{1k}A_k(t,x) v(t,x) + \tilde{q}(t,x)v(t,x)\right|^2 dV_{g}dt \nonumber\\
				&\leq 8\|A\|^2_{L^\infty(M_T)}\|\nabla_g v\|^2_{L^2(M_T)} + 8\lambda^2\|A\|^2_{L^\infty(M_T)}\| v\|^2_{L^2(M_T)} + 2\|\tilde{q}\|^2_{L^\infty(M_T)}\| v\|_{L^2(M_T)}^2
			\end{align}
			Combining  $I_2$, $I_4$, $I_5$,  $I_6$, and $I_{s, 2}$ in the following way:
			\begin{align*}
				& I_2 + I_4 + I_5+ I_6 - I_{s, 2}  =\lambda^2\left(s\left(3 \beta^2 - 1\right) - 8\|A\|^2_{L^\infty(M_T)} - \frac{2}{\lambda^2}\|\tilde{q}\|^2_{L^\infty(M_T)}\right)\|v\|_{L^2(M_T)}^2 \\ & \qquad  + \lb  \mbox{terms having lower power of } \lambda \rb \|v\|_{L^2(M_T)}^2  + \left(3s -  8\|A\|^2_{L^\infty(M_T)}\right) \|\nabla_g v\|^2_{L^2(M_T)} \\
				&\qquad +\int_{\Sigma} (\lambda-s (x_1 + 2 \ell))\lvert \partial_{\nu}v\rvert_{g}^{2}\langle \nu,e_{1}\rangle_{g}\ dS_{g}dt.
			\end{align*}
			After choosing $s$ and $\lambda$ large enough together with using the fact that $\beta \in (1/\sqrt{3}, 1) $, and a combination of  all estimates obtained above, will amount to have the  following estimate on $\left\|P_s v\right\|_{L^2(M_T)}^2 $:
			\begin{align}\label{Estimate for conjugated operator}
				\begin{aligned}
					\left\|P_s v\right\|_{L^2(M_T)}^2 &\geq C\left(s\lambda^2  \|v\|_{L^2(M_T)}^2 + \|\nabla_g v(T, \cdot)\|_{L^2(M)}^2 -\lambda^2  \|v(T, \cdot )\|_{L^2(M)}^2 
					+ s  \left\|  \nabla_g v\right\|_{L^2(M_T)}^2\right. \\ &\qquad \left.+ \lambda \int_{\Sigma} \lvert \partial_{\nu}v(t,x)\rvert_{g}^{2}\langle \nu,e_{1}\rangle_{g}\ dS_{g}dt\right)
				\end{aligned}
			\end{align}
			where constant $C$ depends only on  $A$, $q$, $T$ and $M$.
			
			This provides the estimate for the operator $P_s = e^{-\varphi_s(t, x)}\mathcal{L}_{A,q} e^{\varphi_s(t, x)}$. To obtain the required Carleman estimate, we put $v(t, x) = e^{-\varphi_s(t, x)} u(t, x)$ 
			\begin{align*}
				&\| e^{-\varphi_s} \mathcal{L}_{A,q} u\|_{L^2(M_T)}^2 +\lambda^2 C \|e^{-\varphi_s(T, \cdot)}u(T, \cdot )\|_{L^2(M)}^2 \ +\lambda \int_{ \Sigma_{-}} \lvert e^{-\varphi_s} \partial_{\nu}u(t,x)\rvert_{g}^{2}|\langle \nu,e_{1}\rangle_{g}| \ dS_{g}dt\\
				& \quad \geq s\lambda^2  \|e^{-\varphi_s} u\|_{L^2(M_T)}^2 + \|e^{-\varphi_s(T, \cdot)} \nabla_g u(T, \cdot)\|_{L^2(M)}^2  + s  \left\|e^{-\varphi_s} \nabla_g u\right\|_{L^2(M_T)}^2\\
				&\qquad +\lambda\int_{ \Sigma_{+}} \lvert e^{-\varphi_s} \partial_{\nu}u(t,x)\rvert_{g}^{2}\langle \nu,e_{1}\rangle_{g}\ dS_{g}dt.
			\end{align*}
			Finally, using the expression for $\varphi(t,x)$ and the  fact that $e^{-\frac{s\lb x_{1}+2\ell\rb^{2}}{2}}$ has a strictly positive lower and upper bound, we get the following required estimate
			\begin{align*}
				&\lambda^2  \|e^{-\varphi} u\|_{L^2(M_T)}^2 +\left\|e^{-\varphi} \nabla_g u\right\|_{L^2(M_T)}^2+ \|e^{-\varphi(T, \cdot)} \nabla_g u(T, \cdot)\|_{L^2(M)}^2   +\lambda\int_{ \Sigma_{+}} e^{-2\varphi}\lvert  \partial_{\nu}u(t,x)\rvert^{2}\langle \nu,e_{1}\rangle_{g}\ dS_{g}dt\\
				&\ \ \leq C\lb \| e^{-\varphi} \mathcal{L}_{A,q} u\|_{L^2(M_T)}^2 +\lambda^2  \|e^{-\varphi(T, \cdot)}u(T, \cdot )\|_{L^2(M)}^2 +\lambda \int_{ \Sigma_{-}} e^{-2\varphi}\lvert  \partial_{\nu}u(t,x)\rvert^{2}| \langle \nu,e_{1}\rangle_{g}| \ dS_{g}dt\rb
			\end{align*}
			for some constant $C>0$ independent of $\lambda$. This completes the proof of the theorem. 
		\end{proof}
	\end{theorem}
	Our next aim of this section is to derive the interior Carleman estimates in a Sobolev space of negative order for $\mathcal{L}_{A,q}$ and its formal  $L^{2}$-adjoint $\mathcal{L}_{A,q}^{*}$.  
	Before going to state and prove the interior Carleman estimates, we first give some definitions and notations for semi-classical Sobolev spaces of arbitrary order. This will help us to represent the Carleman estimates in a nice form. 
	Let us begin by  assuming that $(M,g)$ is embedded in a compact Riemannian manifold $(N,g)$ without boundary and denote by $N_{T}:=(0,T)\times N.$ Following \cite{Ferreira_Kenig_Salo_Uhlmann_invention}, we denote by $J^{s}$ for $s\in \mathbb{R}$,  the semi-classical pseudo-differential operator  of order $s$ on $(N,g)$  and it is defined by  $\displaystyle J^{s}:=\left(\lambda^{2}-\Delta_{g}\right)^{s/2}$. Using this we define the semi-classical Sobolev space $H^{s}_{\lambda}(N)$ for $s\in \mathbb{R}$, as the completion of $C^{\infty}(N)$ with respect to the following norm
	\[\lVert u\rVert_{H^{s}_{\lambda}(N)}:=\lVert J^{s}u\rVert_{L^{2}(N)}.\] Since $(N,g)$ is a Riemannian manifold without boundary therefore the dual of $H^{s}_{\lambda}(N)$,  for any $s\in\mathbb{R}$ can be identified with $H^{-s}_{\lambda}(N)$. Also note that for $s=1$, we have that \[\lVert u\rVert^{2}_{H^{1}_{\lambda}(N)}:=\lambda^{2} \lVert u\rVert^{2}_{L^{2}(N)}+\lVert \nabla_{g}u\rVert^{2}_{L^{2}(N)}.\]
	Now following \cite{Evans_Book} the time-dependent Sobolev sapace $L^{2}(0,T;H^{s}_{\lambda}(N))$ is defined as the set of all strongly measurable functions  $u:[0,T]\rightarrow H^{s}_{\lambda}(N)$ such that 
	\begin{align}\label{Norm on tieme-dependent Sobolev}
		\lVert u\rVert_{L^{2}(0,T;H^{s}_{\lambda}(N)}:=\left(\int_{0}^{T}\lVert u(t,\cdot)\rVert_{H^{s}_{\lambda}(N)}^{2}\ dt\right)^{1/2}<\infty.
	\end{align} 
	Then $L^{2}(0,T;H^{s}_{\lambda}(N))$ is a Banach space with respect to the norm $\lVert\cdot\rVert_{L^{2}(0,T;H^{s}_{\lambda}(N))}$ defined by \eqref{Norm on tieme-dependent Sobolev} and the dual of $L^{2}(0,T;H^{s}_{\lambda}(N))$ can be identified with $L^{2}(0,T;H^{-s}_{\lambda}(N))$. 
	If we take  $v\in C_{c}^{\infty}(M_{T})$ in \eqref{Estimate for conjugated operator} then we have the following estimate 
	\begin{align}\label{H1-L2 estimate}
		\begin{aligned}
			\lVert v\rVert^{2}_{L^{2}(0,T;H^{1}_{\lambda}(N))}\leq C\lVert \mathcal{L}_{\varphi}v\rVert^{2}_{L^{2}(0,T;L^{2}(N))}
		\end{aligned}
	\end{align}
	where  $\varphi$ is same as in Theorem \ref{Boundary Carleman estimate Theorem} and $\displaystyle \mathcal{L}_{\varphi}:= e^{-\varphi}\mathcal{L}_{A,q}e^{\varphi}$. Now if we denote by $\mathcal{L}_{\varphi}^{*}:= e^{\varphi}\mathcal{L}_{A,q}^{*}e^{-\varphi}$  where $\mathcal{L}_{A,q}^{*}$ stands for a formal $L^{2}-$adjoint of $\mathcal{L}_{A,q}$ then using the arguments similar to the one used in deriving \eqref{Estimate for conjugated operator}, the following estimate
	\begin{align}\label{H1-L2 estimate for adjoint}
		\lVert u\rVert^{2}_{L^{2}(0,T;H^{1}_{\lambda}(N))}\leq C\lVert \mathcal{L}^{*}_{\varphi}u\rVert^{2}_{L^{2}(0,T;L^{2}(N))}
	\end{align}
	holds for all $u\in C_{c}^{\infty}(M_{T})$ where $\varphi$ is same as in Theorem \ref{Boundary Carleman estimate Theorem} and constant $C>0$ is independent of $\lambda$ and $u$. 
	
	In order to construct the suitable solutions to $\mathcal{L}^{*}_{A,q}u=0$ and $\mathcal{L}_{A,q}v=0$, we need to shift the index by $-1$ for spacial variable in  \eqref{H1-L2 estimate} and \eqref{H1-L2 estimate for adjoint} respectively, which we will do in the following lemma. 
	\begin{lemma}\label{Interior estimate in negative order}
		Let $\displaystyle \mathcal{L}_{\varphi}^{*}:= e^{\varphi}\mathcal{L}_{A,q}^{*}e^{-\varphi},$ and $\displaystyle \mathcal{L}_{\varphi}:= e^{-\varphi}\mathcal{L}_{A,q}e^{\varphi},$ where $\mathcal{L}_{A,q}^{*}$ denote the formal $L^{2}$-adjoint of $\mathcal{L}_{A,q}$ and $\varphi$, $A$ and $q$ be as in Theorem \ref{Boundary Carleman estimate Theorem}. Then there exists a constant $C>0$ independent of $\lambda$  and $v$ such that \begin{align}\label{Interior estimate for adjoint operator}
			\lVert v\rVert_{L^{2}([0,T];L^{2}(N))}\leq C\lVert \mathcal{L}_{\varphi}^{*}v\rVert_{L^{2}(0,T;H^{-1}_{\lambda}(N))}
		\end{align}
		holds for all $\lambda$ large enough  and for all $v\in C_{c}^{\infty}(M_T)$ and 
		\begin{align}\label{Interior estimate for operator}
			\lVert v\rVert_{L^{2}([0,T];L^{2}(N))}\leq C\lVert \mathcal{L}_{\varphi}v\rVert_{L^{2}(0,T;H^{-1}_{\lambda}(N))}
		\end{align}
		holds for all $\lambda$ large enough and for all $v\in C_{c}^{\infty}(M_T)$. 
		\begin{proof}
  First, we establish \eqref{Interior estimate for adjoint operator}, and the proof for \eqref{Interior estimate for operator} can be carried out in a similar manner.
We begin with the inequality:
   \begin{align*}
		\lVert v\rVert^{2}_{L^{2}(0,T;H^{1}_{\lambda}(N))}\leq C\lVert \mathcal{L}^{*}_{\varphi}v\rVert^{2}_{L^{2}(0,T;L^{2}(N))}
  \end{align*}
holds for all $v\in C_{c}^{\infty}(M_{T})$.
   Next, we shift the index by $-1$ in the above estimate. Let $w \in  C_{c}^{\infty}(M_T)$ and consider the adjoint operator defined as:
 \[\displaystyle\mathcal{L}^{*}_{A,q}:=\Big(-  \partial_{t}-\sum_{j,k=1}^{n}\frac{1}{\sqrt{\lvert g\rvert}}( \partial_{x_{j}}-A_{j})(\sqrt{\lvert g\rvert} g^{jk}( \partial_{x_{k}}- A_{k}))+\overline{q}\Big).\]
 For $s>0$, define the convexified weight function $\varphi_{s}$ as follows:
 \begin{align*} 
				\varphi_s(t, x) := \varphi(t, x) + \frac{s(x_1 + 2 \ell)^2}{2}=\lambda^{2}\beta^{2}t+\lambda x_{1}+ \frac{s(x_1 + 2 \ell)^2}{2}.
		\end{align*} 
   Let $\displaystyle P^{*}_{s}:= e^{\varphi_s}\mathcal{L}_{A,q}^{*}e^{-\varphi_s}$, we have
   \begin{align*}
        P_s^{*}w = e^{\varphi_s}\lb -\partial_{t}-\partial_{x_{1}}^{2}-\Delta_{g_0}+2\langle A,\nabla_{g}\rangle_{g}+\tilde{q}^* \rb (e^{-\varphi_s}w)
  \end{align*}
   where $\tilde{q}^*(t, x) := \overline{q}(t, x) + \delta_g A(t, x) - \lvert A(t,x)\rvert_g^2$.\\
   Expressing $P_s^{*}w$ as a sum of three components:
   \begin{align*}
       P_s^{*}w := P^*_1w(t,x)+P^*_2w(t,x)+P^*_3w(t,x)
   \end{align*}
   where 
   \begin{align*}
   \begin{aligned}
					P_1^{*} w(t,x)&:= \left(-\partial_{t}w + 2(\lambda+ s (x_1 + 2 \ell))\partial_{x_1}w + 4sw\right)(t,x),\\
					P_2^{*} w(t,x)&:= \left (-\partial^2_{x_1}w - \Delta_{g_0} w - \lambda^2(1 - \beta^2)w - 2 \lambda s(x_1 + 2 \ell) w - s^2 (x_1 + 2 \ell)^2 w -  3sw\right)(t,x),\\
					P_3^{*} w(t,x)&:= 2 \langle A(t,x) , \nabla_g w(t,x)\rangle_{g} -2 \lb \lambda + s(x_1 + 2\ell) \rb g^{1k}A_k(t,x) w(t,x) + \tilde{q}^{*}(t,x)w(t,x).
				\end{aligned}
			\end{align*}
Furthermore, we have $\displaystyle J^{-1}(P^*_1+P^*_2)J^1 w = (P^*_1+P^*_2)w$, from this we get
   \begin{align*}
       \|(P^*_1+P^*_2)J^1 w\|_{L^2(0,T;H_{\lambda}^{-1}(N))} &= \|J^{-1}(P^*_1+P^*_2)J^1 w\|_{L^2(0,T;L^2(N))}\\
       &= \|(P^*_1+P^*_2) w\|_{L^2(0,T;L^2(N))}.
   \end{align*}
   From the same calculation as done for the Carleman estimate \eqref{eq: BC estimate}, we obtain
   \begin{align*}
   \|(P^*_1+P^*_2)J^1 w\|_{L^2(0,T;H_{\lambda}^{-1}(N))}\geq \sqrt{s}\|\nabla_{g}w\|_{L^2(0,T;L^2(N))} + \sqrt{s}\lambda \|w\|_{L^2(0,T;L^2(N))}. 
   \end{align*}
   Now consider
    \begin{align*}
   \|P^*_3 J^1 w\|_{L^2(0,T;H_{\lambda}^{-1}(N))}\leq C(\|A\|_{\infty}\|\nabla_{g}w\|_{L^2(0,T;L^2(N))} + \lambda \|A\|_{\infty}\|w\|_{L^2(0,T;L^2(N))} + \|\tilde{q}^*\|_{\infty}\|w\|_{L^2(0,T;L^2(N))}) 
   \end{align*}
Hence using the inequality and choosing $s$ and $\lambda$ large enough, we get
   \begin{align*}
   \| P_s^{*}J^1 w\|_{L^2(0,T;H_{\lambda}^{-1}(N))}\geq C \|w\|_{L^2(0,T;H_{\lambda}^{1}(N))}.  
   \end{align*}
Now, consider $\chi \in C^{\infty}_c(\widetilde{M})$ such that $\chi = 1$ in $\overline{M_1}$ where $\overline{M}\subset M_1\subset\widetilde{M}$. By taking $w = \chi J^{-1} v$ in the above estimate and using
\begin{align*}
    \|(1-\chi)J^{-1}v\|_{L^2(0,T;H^{1}_{\lambda}(N))} \leq \frac{C}{\lambda^2}\|v\|_{L^2(0,T;L^2(N))}
\end{align*}
and 
\begin{align*}
    \|v\|_{L^2(0,T;L^2(N))} &= \|J^{-1}v\|_{L^2(0,T;H^{1}_{\lambda}(N))}\\
    &\leq \|w\|_{L^2(0,T;H^{1}_{\lambda}(N))} + \frac{C}{\lambda^2}\|v\|_{L^2(0,T;L^2(N))}.
\end{align*}
We get
\begin{align*}
    \| P_s^{*} v\|_{L^2(0,T;H^{-1}_{\lambda}(N))} &\geq \| P_s^{*} J^1 w\|_{L^2(0,T;H^{-1}_{\lambda}(N))} -  \frac{C}{\lambda^2}\|v\|_{L^2(0,T;L^2(N))}\\
    & \geq \|w\|_{L^2(0,T;H^{1}_{\lambda}(N))}  -  \frac{C}{\lambda^2}\|v\|_{L^2(0,T;L^2(N))}\\
    & \geq C \|v\|_{L^2(0,T;L^2(N))}
\end{align*}
hold for $\lambda$ large. Therefore, we have
\begin{align*}
		\lVert v\rVert_{L^{2}([0,T];L^{2}(N))}\leq C\lVert P_s^{*}v\rVert_{L^{2}(0,T;H^{-1}_{\lambda}(N))}.
\end{align*}
Now, using the expression for $\varphi(t,x)$ and the fact that $e^{\frac{s\lb x_{1}+2\ell\rb^{2}}{2}}$ has a strictly positive lower and upper bound, we conclude
\begin{align*}
		\lVert v\rVert_{L^{2}([0,T];L^{2}(N))}\leq C\lVert \mathcal{L}_{\varphi}^{*}v\rVert_{L^{2}(0,T;H^{-1}_{\lambda}(N))}
\end{align*}
		holds for all $\lambda$ large enough  and for all $v\in C_{c}^{\infty}(M_T)$.
  \end{proof}
	\end{lemma}
	The above estimates, together with the Hahn-Banach theorem and the Riesz representation theorem, give the following solvability result, proof of which follows from \cite{Suman_Manmohan_IPI,Soumen_Manmohan_EECT}.  
	\begin{lemma}\label{Existence of solution lemma}
		Let $\varphi$, $A$ and $q$ be as before and $\lambda>0$ be large enough. Then for $F\in L^{2}(M_{T})$ there exists a solution $u\in H^{1}(0,T;H^{-1}(M))\cap L^{2}(0,T;H^{1}(M))$  of 
		\begin{align*}
			\begin{aligned}
					\mathcal{L}_{\varphi}w(t,x)=F(t,x),\  (t,x)\in M_{T}
			\end{aligned}
		\end{align*}
		satisfying the following estimate
		\begin{align}
			\lVert u\rVert_{L^{2}(0,T;H^{1}_{\lambda}(M))}\leq C\lVert F\rVert_{L^{2}(M_{T})}  
		\end{align}
		for some constant $C>0$ independent of $\lambda$ and $u$ and there exists a solution $v\in H^{1}(0,T;H^{-1}(M)) \cap L^{2}(0,T;H^{1}(M))$  of 
		\begin{align*}
			\begin{aligned}
					\mathcal{L}_{\varphi}^{*}w(t,x)=F(t,x),\  (t,x)\in M_{T}
			\end{aligned}
		\end{align*}
		satisfying the following estimate
		\begin{align}
			\lVert v\rVert_{L^{2}(0,T;H^{1}_{\lambda}(M))}\leq C\lVert F\rVert_{L^{2}(M_{T})}  
		\end{align}
		for some constant $C>0$ independent of $\lambda$ and $v$.  
  \begin{proof}
The proof for $\mathcal{L}_{\varphi}$ is presented below, and the proof for $\mathcal{L}_{\varphi}^{*}$ can be established using analogous arguments.\\
Consider the subspace $S$ of ${L^{2}(0,T;H^{-1}_{\lambda}(N))}$ defined as
\begin{align*}
S:=\{\mathcal{L}_{\varphi}^{*}w(t,x): w\in C_{c}^{\infty}(M_T) \}.
\end{align*}
Define the linear operator $T$ on $S$ by
\begin{align*}
 T(\mathcal{L}_{\varphi}^{*}z) = \int_{M_T} z(t,x)\overline{F(t,x)}\ dV_g dt, \ \text{for} \ F\in L^2(M_T).
\end{align*}
For any $\mathcal{L}_{\varphi}^{*}z \in S$, we have
\begin{align*}
|T(\mathcal{L}_{\varphi}^{*}z)| \leq \int_{M_T} |z(t,x)| |F(t,x)| \ dV_g dt  \leq \|z\|_{L^2(M_T)} \|F\|_{L^2(M_T)}.
\end{align*}
Using the Carleman estimate \eqref{Interior estimate for adjoint operator}, we obtain
\begin{align*}
|T(\mathcal{L}_{\varphi}^{*}z)| \leq C \|F\|_{L^2(M_T)} \| \mathcal{L}_{\varphi}^{*}z\|_{L^{2}(0,T;H^{-1}_{\lambda}(N))}.
\end{align*}
This inequality holds for $z\in C^{\infty}_c(M_T)$. By the Hahn-Banach theorem, extend the linear operator $T$ to $L^{2}(0, T; H^{-1}_{\lambda}(N))$. Denote the extended map as $T$ and note that it satisfies the inequality
\begin{align*}
  \|T\| \leq C \|F\|_{L^2(M_T)}.
\end{align*}
By the Riesz representation theorem, as $T$ is a bounded linear functional on $L^{2}(0, T; H^{-1}_{\lambda}(N))$, there exists a unique $u \in L^{2}(0, T; H^{1}_{\lambda}(N))$ such that
\begin{align*}
T(f) = \langle f, u\rangle_{L^{2}(0,T;H^{-1}_{\lambda}(N)),L^{2}(0,T;H^{1}_{\lambda}(N))} \ \text{for} \ f\in L^{2}(0,T;H^{-1}_{\lambda}(N)),
\end{align*}
with
\begin{align*}
 \|u\|_{L^{2}(0,T;H^{1}_{\lambda}(N))} \leq C\|F\|_{L^{2}(M_T)}.
\end{align*}
Now, for $z\in C^{\infty}_c(M_T)$, choosing $f = \mathcal{L}_{\varphi}^{*}z$ in the above equation, we get $\mathcal{L}_{\varphi}u = F$.\\
Using the expression for $\mathcal{L}_{\varphi}$ and the fact that $u\in L^{2}(0,T;H^{1}(M))$ and $F\in L^{2}(M_T) $, we conclude that $\partial_t u \in L^{2}(0,T;H^{-1}(M))$. Hence, we have $ u\in H^{1}(0,T;H^{-1}(M))\cap L^{2}(0,T;H^{1}(M))$.\\
 \end{proof}
\end{lemma}

\section{Construction of geometric optics solutions}\label{Construction of GO solutions}
	In this section, we aim to construct the exponential growing and decaying solutions to the convection-diffusion operator $\displaystyle \mathcal{L}_{A,q}$ and its $L^{2}$-adjoint $\displaystyle \mathcal{L}^{*}_{A,q}$, respectively. Construction of these solutions will be proved with the help of the interior  Carleman estimates in negative order Sobolev spaces stated in Lemma \ref{Interior estimate in negative order}.  
	\subsection{Construction of exponentially growing solutions}
	In this subsection, we will construct the exponential growing solutions to  $\mathcal{L}_{A,q}u(t,x)=0,$ in $M_{T}$ which takes  the following form
	\begin{align}\label{Growing expression}
		u(t,x)=e^{\lb \varphi+i\psi\rb(t,x)}\Big (T_{g}(t,x)+R_{g,\lambda}(t,x)\Big)
	\end{align}
	where $\varphi$ is same as  in Theorem \ref{Boundary Carleman estimate Theorem} and  $\psi$,  $T_{g}$ will be constructed using the WKB construction in such a way that  the correction term  $R_{g,\lambda}$  satisfies the following
	\begin{align*}
		\begin{aligned}
			e^{-\lb \varphi+i\psi\rb}\mathcal{L}_{A,q}\left(e^{\lb \varphi+i\psi\rb}R_{g,\lambda}(t,x)\right)&=F_{\lambda}(t,x),\ (t,x)\in M_{T}\\
		\end{aligned}
	\end{align*}
	for $F_{\lambda}\in L^{2}(M_{T})$ such that $\lVert F_\lambda\rVert_{L^{2}(M_{T})}\leq C,$ for some constant $C>0$ independent of $\lambda$ and $R_{g,\lambda}$ satisfies 
	$\displaystyle \lVert R_{g,\lambda}\rVert_{L^{2}(0,T;H^{1}_{\lambda}(M))}\leq C \lVert F_{\lambda}\rVert_{L^{2}(M_{T})},$ for some constant $C>0,$ not depending on $\lambda$. More precisely, we prove the following theorem. 
\begin{theorem}\label{Growing solutions Theorem}
Let $M_{T}$, $\mathcal{L}_{A,q}$ and $\varphi$ be as before. Suppose $(D,g_{0})$ be a simple  manifold which is extension of $(M_{0},g_{0})$ in the sense that $M_{0}\subset D$ and $y_{0}\in  D$ is such that $(x_{1},y_{0})\notin M$ for all $x_{1}$. Now if $(r,\theta)$ denote the polar normal coordinates on $(D,g_{0})$, $(x_{1},r,\theta)$ denote the points in $M$ and $A_{1}$ and $A_{r}$ are components of $A$ in $x_{1}$ and $r$ coordinates respectively,  then for $\lambda$ large enough the following  equation  
\begin{align*}
\begin{aligned}
\mathcal{L}_{A,q}v(t,x)=0,\ \ (t,x)\in M_{T}
\end{aligned}
\end{align*}
has a solution taking the following form
\begin{align}\label{eq:u in polar coordinates}
\begin{aligned}
u(t,x)=e^{\varphi+i\psi}\Big (T_{g}(t,x_{1},r,\theta)+R_{g,\lambda}(t,x_{1},r,\theta)\Big)
\end{aligned}
\end{align}
where \begin{align*}
\begin{aligned}
\psi=\lambda(\sqrt{1-\beta^{2}})r,\ \mbox{and} \  T_{g}(t,x_{1},r,\theta)=\phi(t)e^{i\mu\left(\sqrt{1-\beta^{2}}\right)x_{1}}e^{-\mu r}e^{i\Phi_{1}(t, x_{1},r,\theta)}b(r,\theta)^{-1/4}h(\theta) 
\end{aligned}
\end{align*}  
here $\phi\in C_{c}^{\infty}(0,T)$, $\Phi_{1}$ is solution to  $$\partial_{1}\Phi_{1}+i(\sqrt{1-\beta^{2}})\partial_{r}\Phi_{1}+\lb - i A_{1}+(\sqrt{1-\beta^{2}})A_{r}\rb=0$$ and $R_{g,\lambda}$ satisfies the following 
\begin{align*}
			\begin{aligned}
					\mathcal{L}_{\varphi}\lb e^{i\psi}R_{g,\lambda}\rb(t,x)=-e^{i\psi}\mathcal{L}_{A,q}T_{g}(t,x),\ \ \ (t,x)\in M_{T}
			\end{aligned}
		\end{align*} and $\displaystyle \lVert R_{g,\lambda}\rVert_{L^{2}\lb 0,T;H^{1}_{\lambda}(M)\rb}\leq C$ for some constant $C>0$ independent of $\lambda$. 
		\begin{proof}
			Following \cite{Ferreira_Kenig_Salo_Uhlmann_invention}, if we denote  $\rho:=\varphi+i\psi,$ then simple calculations show that  the conjugated
			operator $\mathcal{L}_{\rho}:=e^{-\rho}\mathcal{L}_{A,q}e^{\rho}$ will have the following expression
			\begin{align*}
				\begin{aligned}
					\mathcal{L}_{\rho}=\mathcal{L}_{A,q}+\left(\partial_{t}\rho-\Delta_{g}\rho-g^{jk}\partial_{j}\rho\partial_{k}\rho\right)-2\left(g^{jk}\partial_{j}\rho\partial_{k}+g^{jk}\partial_{j}\rho A_{k}\right). 
				\end{aligned}
			\end{align*}
			Using $\rho=\varphi+i\psi,$ and  $\varphi=\lambda^{2}\beta^{2}t+\lambda x_{1}$, we get
			\begin{align}\label{Expression for conjugated operator with rho}
				\begin{aligned}
	    \mathcal{L}_{\rho}&=\mathcal{L}_{A,q}+\left(\lambda^{2}\beta^{2}-\lambda^{2}+g^{jk}\partial_{j}\psi\partial_{k}\psi\right)\\
					&\ \ \ \ \ \ \ \ \ \  -\left(2\lambda \partial_{1}+2ig^{jk}\partial_{j}\psi\partial_{k}+2\lambda A_{1}+2 i g^{jk}\partial_{j}\psi A_{k}+i\Delta_{g}\psi + 2i\lambda \partial_1\psi-i\partial_t\psi\right).
				\end{aligned}
			\end{align}
			Now $u$ given by \eqref{Growing expression} solves $\mathcal{L}_{A,q}v=0$ if and only if   $\mathcal{L}_{\rho}\lb e^{-\rho}u\rb=0$. This will give us 
			\begin{align}\label{Equation for Rdlambda growing}
				\begin{aligned}
					\mathcal{L}_{\rho}R_{g,\lambda}(t,x)&=-\mathcal{L}_{A,q}T_{g}(t,x)-\left(\lambda^{2}\beta^{2}-\lambda^{2}+g^{jk}\partial_{j}\psi\partial_{k}\psi\right)T_{g}(t,x)\\
					&\ \ \ \ +\left(2\lambda \partial_{1}+2ig^{jk}\partial_{j}\psi\partial_{k}+2\lambda A_{1}+2ig^{jk}\partial_{j}\psi A_{k}+i\Delta_{g}\psi {+ 2i\lambda \partial_1\psi-i\partial_t\psi}\right)T_{g}(t,x),\ \ (t,x)\in M_{T}. 
				\end{aligned}
			\end{align}
			In order to have $\lVert R_{g,\lambda}\rVert_{L^{2}(0,T;H^{1}_{\lambda}(M)}\leq C$, we choose $\psi$ and $T_{g}$ such that 
			\begin{align}\label{Eikonal equation}
                 {\partial_t\psi  =  0}, \ \ g^{jk}\partial_{j}\psi\partial_{k}\psi=\lambda^{2}(1-\beta^{2})
			\end{align}
			and 
			\begin{align}\label{Transport equation}
				\left(2\lambda \partial_{1}+2ig^{jk}\partial_{j}\psi\partial_{k}+2\lambda A_{1}+2ig^{jk}\partial_{j}\psi A_{k}+i\Delta_{g}\psi{+ 2i\lambda\partial_1\psi}\right)T_{g}(t,x)=0,\ (t,x)\in M_{T}. 
			\end{align}
			To solve equations \eqref{Eikonal equation} and \eqref{Transport equation} for $\psi$ and $T_{g}$, we use the polar normal coordinates $(r,\theta)$ on $(D,g_{0})$ centered at $y_{0}\in D$ as mentioned in statement of  theorem. 
			We consider the polar normal coordinates on $D$ which are  denoted  by $(r,\theta)$ and given by $x_{0}=\exp_{y_{0}}(r\theta),$ where $r>0$ and $\theta\in S_{y_{0}}(D):=\{v\in T_{y_0}D:\ \lvert v\rvert_{g}=1\}$, here $T_{y_{0}}D$ denote the tangent space to $D$ at $y_{0}\in D$. Then using the Gauss lemma (see Lemma $15$ in Chapter $9$ of \cite{Spivak_DG_Vol 1}) there exists a smooth positive definite matrix $P(r,\theta)$ with $\det(P):=b(r,\theta)$ such that   the metric $g_{0}$ in the polar normal coordinates $(r,\theta)$, takes the following form
			\begin{align}\label{Form of g0 in polar normal}
				\begin{aligned}
					g_{0}(r,\theta)=\begin{bmatrix}
						1&0\\
						0&P(r,\theta)
					\end{bmatrix}.
				\end{aligned}
			\end{align}
			Now since the points in $M$ are denoted by $(x_{1},r,\theta)$ where $(r,\theta)$ are polar normal coordinates in $(D,g_{0})$,  therefore  after using the previous form of $g_{0}$, the metric $g$ has the following form
			\begin{align}\label{form of g after polar normal coordinates}
				\begin{aligned}
					g(x_{1},r,\theta)=\begin{bmatrix}
						1 & 0 & 0\\
						0 & 1 & 0\\
						0 & 0 & P(r,\theta)
					\end{bmatrix}.
				\end{aligned}
			\end{align}
			Using \eqref{form of g after polar normal coordinates}, we see that \begin{align}\label{Solution to Eikonal equation}
				\psi(x)=(\lambda\sqrt{1-\beta^{2}})\text{dist}_{g}(y_{0},x)=(\lambda\sqrt{1-\beta^{2}})r,
			\end{align} solves  equation \eqref{Eikonal equation} on $M$. Using this choice of $\psi$ and form of $g$ given by \eqref{form of g after polar normal coordinates} in equation \eqref{Transport equation}, we have 
			\begin{align*}
				\begin{aligned}
					\left( \partial_{1}+i(\sqrt{1-\beta^{2}})\partial_{r}+ A_{1}+i(\sqrt{1-\beta^{2}})A_{r}+i(\sqrt{1-\beta^{2}})\frac{\partial_{r}     b(r,\theta)}{4b(r,\theta)}\right)T_{g}(t,x_1,r,\theta)=0.
				\end{aligned}
			\end{align*}
			Now, one can check that the solution of the above equation can be given by 
			\begin{align}\label{Expression for Tg}
				T_{g}(t,x_{1},r,\theta)=\phi(t)e^{i\mu \left(\sqrt{1-\beta^{2}}\right)x_{1}}e^{-\mu r}e^{i\Phi_{1}(t, x_{1},r,\theta)}b(r,\theta)^{-1/4}h(\theta)  
			\end{align}
			where $\phi\in C_{c}^{\infty}(0,T)$, $\mu\in\rr$, $h\in C^{\infty}(S_{y_{0}}(D))$ are arbitrary but fixed and $\Phi_{1}(t, x_{1},r,\theta)$ satisfies the  following
			\begin{align}\label{Equation for Phi}
				\lb \partial_{1}\Phi_{1}+i\lb\sqrt{1-\beta^{2}}\rb\partial_{r}\Phi_{1}\rb +\lb -i A_{1}+(\sqrt{1-\beta^{2}})A_{r}\rb=0. 
			\end{align}
			Now using \eqref{Solution to Eikonal equation} and \eqref{Expression for Tg} in \eqref{Equation for Rdlambda growing}, we get 
			\begin{align*}
				\begin{aligned}
						\mathcal{L}_{\rho}R_{g,\lambda}(t,x)=-\mathcal{L}_{A,q}T_{g}(t,x),\ \ (t,x)\in M_{T}
				\end{aligned}
			\end{align*}
			But $\mathcal{L}_{\rho}R_{g,\lambda}=e^{-i\psi}\mathcal{L}_{\varphi}\left(e^{i\psi}R_{g,\lambda}\right)$ therefore if we denote $\widetilde{R}_{g,\lambda}=e^{i\psi}R_{g,\lambda}$ then $\widetilde{R}_{g,\lambda}$ satisfies the following equation
			\begin{align}\label{Equation for tilde Rglambda}
					\mathcal{L}_{\varphi}\widetilde{R}_{g,\lambda}(t,x)=-e^{i\psi}\mathcal{L}_{A,q}T_{g}(t,x),\ \ (t,x)\in M_{T}
			\end{align}
			Now using the expressions for $\psi$ and $T_{g}$ from \eqref{Solution to Eikonal equation} and \eqref{Expression for Tg} respectively and assumptions on $A$ and $q$, we have that   $\displaystyle -e^{i\psi}\mathcal{L}_{A,q}T_{g}\in L^{2}(M_{T})$ and $\displaystyle \lVert e^{i\psi}\mathcal{L}_{A,q}T_{g}\rVert_{L^{2}(M_{T})}\leq C$, for some constant $C>0$ independent of $\lambda$. Hence using Lemma \ref{Existence of solution lemma} together with above estimate for right hand side of \eqref{Equation for tilde Rglambda}, we conclude that there exists $\widetilde{R}_{g,\lambda}\in H^{1}(0,T;H^{-1}(M))\cap L^{2}(0,T;H^{1}(M))$ solving \eqref{Equation for tilde Rglambda} and it satisfies the following estimate $\displaystyle \lVert \widetilde{R}_{g,\lambda}\rVert_{L^{2}(0,T;H^{1}_{\lambda}(M))}\leq C$, for some constant $C>0,$ independent of $\lambda$. Hence, we conclude that $R_{g,\lambda}$ solves the required equation and satisfies the desired estimate. This completes the proof of the Theorem. 
		\end{proof}
	\end{theorem}
	\subsection{Construction of exponentially decaying solutions} The aim of this subsection is to construct the exponential decaying solutions to  \[\displaystyle\mathcal{L}^{*}_{A,q}u:=\Big(-  \partial_{t}-\sum_{j,k=1}^{n}\frac{1}{\sqrt{\lvert g\rvert}}( \partial_{x_{j}}-A_{j})(\sqrt{\lvert g\rvert} g^{jk}( \partial_{x_{k}}- A_{k}))+\overline{q}\Big) u=0, \ \mbox{in}\  M_{T}\] taking the following form
	\begin{align}\label{decaying expression}
		u(t,x)=e^{-\lb \varphi-i\psi\rb(t,x)}\Big (T_{d}(t,x)+R_{d,\lambda}(t,x)\Big)
	\end{align}
	where $\varphi$ is the same as  in Theorem \ref{Boundary Carleman estimate Theorem} and  $\psi$,  $T_{d}$ will be constructed using the WKB construction in such a way that  the correction term  $R_{d,\lambda}$  satisfies the following
	\begin{align*}
		\begin{aligned}
				e^{\lb \varphi-i\psi\rb}\mathcal{L}^{*}_{A,q}\left(e^{-\lb \varphi-i\psi\rb}R_{d,\lambda}(t,x)\right)=F_{\lambda}(t,x),\ (t,x)\in M_{T}
		\end{aligned}
	\end{align*}
	for some $F_{\lambda}\in L^{2}(M_{T})$ such that $\lVert F_\lambda\rVert_{L^{2}(M_{T})}\leq C,$ for some constant $C>0$ independent of $\lambda$ and $R_{d,\lambda}$ satisfies 
	$\displaystyle \lVert R_{d,\lambda}\rVert_{L^{2}(0,T;H^{1}_{\lambda}(M))}\leq C \lVert F_{\lambda}\rVert_{L^{2}(M_{T})},$ for some constant $C>0$ not depending on $\lambda$. To construct these solutions, we first start with the construction of $\psi$ and $T_{d}$ following the arguments used in Theorem \ref{Growing solutions Theorem}. Denote by $\rho:=\varphi-i\psi,$ then one can check that the conjugated operator \[ \mathcal{L}_{\rho}^{*}:=e^{\rho}\mathcal{L}_{A,{q}}^{*}e^{-\rho}=e^{\rho}\Big(-  \partial_{t}-\sum_{j,k=1}^{n}\frac{1}{\sqrt{\lvert g\rvert}}( \partial_{x_{j}}-A_{j})(\sqrt{\lvert g\rvert} g^{jk}( \partial_{x_{k}}- A_{k}))+\overline{q}\Big)e^{-\rho}\]is given by 
	\begin{align*}
		\begin{aligned}
			\mathcal{L}_{\rho}^{*}=\mathcal{L}_{A,{q}}^{*}+\lb \partial_{t}\rho-g^{jk}\partial_{j}\rho\partial_{k}\rho\rb+\lb 2g^{jk}\partial_{j}\rho\partial_{k}-2g^{jk}\partial_{j}\rho A_{k}+\Delta_{g}\rho\rb. 
		\end{aligned}
	\end{align*}
	Using $\rho=\varphi-i\psi$ and $\varphi=\lambda^{2}\beta^{2}t+\lambda x_{1}$, we have 
	\begin{align*}
		\mathcal{L}_{\rho}^{*}=\mathcal{L}_{A,{q}}^{*}+\lb \lambda^{2}\beta^{2}-\lambda^{2}+g^{jk}\partial_{j}\psi\partial_{k}\psi\rb +\lb 2\lambda \partial_{1}-2ig^{jk}\partial_{j}\psi\partial_{k}-2\lambda A_{1}+2ig^{jk}\partial_{j}\psi A_{k}-i\Delta_{g}\psi {+ 2i\lambda \partial_1\psi-i\partial_t\psi}\rb. 
	\end{align*}
	Now we observe that  $u$ given by \eqref{decaying expression} solves $\mathcal{L}_{A,q}^{*}v=0$ in $M_{T}$  if and only if $\mathcal{L}^{*}_{\rho}\lb e^{\rho}u\rb=0$ in $M_{T}.$ Using this, we see that $R_{d,\lambda}$ satisfies the following equation
	\begin{align}\label{Equation for Rdlambda}
		\begin{aligned}
			\mathcal{L}_{\rho}^{*}R_{d,\lambda}(t,x)&=-\mathcal{L}_{A,{q}}^{*}T_{d}(t,x)-\lb \lambda^{2}\beta^{2}-\lambda^{2}+g^{jk}\partial_{j}\psi\partial_{k}\psi\rb T_{d}(t,x)\\
			&\ \  \ \ \ -\lb 2\lambda \partial_{1}-2ig^{jk}\partial_{j}\psi\partial_{k}-2\lambda A_{1}+2ig^{jk}\partial_{j}\psi A_{k}-i\Delta_{g}\psi {+ 2i\lambda \partial_1\psi-i\partial_t\psi}\rb T_{d}(t,x). 
		\end{aligned}
	\end{align}
	To get  the estimate $\displaystyle \lVert R_{d,\lambda}\rVert_{L^{2}(0,T;H^{1}_{\lambda}(M))}\leq C,$ for some constant $C>0$ independent of $\lambda$, we choose $\psi$ and $T_{d}$ satisfying the following equations
	\begin{align}\label{Eikonal decaying case}
{\partial_t\psi  =  0}, \ \
		\lambda^{2}\beta^{2}-\lambda^{2}+g^{jk}\partial_{j}\psi\partial_{k}\psi=0
	\end{align}
	and
	\begin{align}\label{Transport decaying case}
		\lb 2\lambda \partial_{1}-2ig^{jk}\partial_{j}\psi\partial_{k}-2\lambda A_{1}+2ig^{jk}\partial_{j}\psi A_{k}-i\Delta_{g}\psi {+ 2i\lambda \partial_1\psi} \rb T_{d}(t,x)=0, \ (t,x)\in M_{T}  
	\end{align}
	respectively. 
	To solve equations \eqref{Eikonal decaying case} and \eqref{Transport decaying case} for $\psi$ and $T_{d}$, we again use the polar normal coordinates $(r,\theta)$ on $(D,g_{0})$ centered at $y_{0}\in  D$ as used in the proof of Theorem \ref{Growing solutions Theorem}. 
	For a fixed $y_{0}\in D$,  we consider the polar normal coordinates on $D$ which are  denoted  by $(r,\theta)$ and given by $x_{0}=\exp_{y_{0}}(r\theta),$ where $r>0$ and $\theta\in S_{y_{0}}(D):\{v\in T_{y_0}D:\ \lvert v\rvert_{g}=1\}$, here $T_{y_{0}}D$ denote the tangent space to $D$ at $y_{0}\in D$. Then using the Gauss lemma (see Lemma $15$ in Chapter $9$ of \cite{Spivak_DG_Vol 1}) there exists a smooth positive definite matrix $P(r,\theta)$ with $\det P(r,\theta)=b(r,\theta)$ such that the metric $g_{0}$ in the polar normal coordinates $(r,\theta)$, takes form given by \eqref{Form of g0 in polar normal}. 
		Now since the points in $M$ are denoted by $(x_{1},r,\theta)$ where $(r,\theta)$ are polar normal coordinates in $(D,g_{0})$,  therefore  after using the form of $g_{0}$ given  by \eqref{Form of g0 in polar normal}, the metric $g$ takes the form given by equation \eqref{form of g after polar normal coordinates} and using this, we observe that 
			\begin{align}\label{Expression for psi decaying case}
				\psi(x)=\lb \lambda\sqrt{1-\beta^{2}}\rb \text{dist}_{g}(y_{0},x)=\lb \lambda\sqrt{1-\beta^{2}}\rb r,
			\end{align}
			solves equation \eqref{Eikonal decaying case} and 
			\begin{align}\label{Expression for Td decaying case}
				T_{d}(t,x_{1},r,\theta)=\phi(t)e^{i\Phi_{2}(t,x_{1},r,\theta)}b(r,\theta)^{-1/4}h(\theta)
			\end{align}
			solves equation \eqref{Transport decaying case} where $\phi\in C_{c}^{\infty}(0,T)$, $h\in C^{\infty}(S_{y_{0}}(D))$ are arbitrary but fixed and $\Phi_{2}(t,x_{1},r,\theta)$ satisfies the following
			\begin{align}\label{Equation for Phi2 decaying case}
				\lb \partial_{1}\Phi_{2}-i\lb \sqrt{1-\beta^{2}}\rb \partial_{r}\Phi_{2}\rb +\lb i A_{1}+ \sqrt{1-\beta^{2}} A_{r}\rb =0, 
			\end{align}
			$A_{1}$ and $A_{r}$ are components of $A$ in $x_{1}$ and $r$ coordinates respectively. Now if we use  \eqref{Expression for psi decaying case} and \eqref{Expression for Td decaying case} in equation \eqref{Equation for Rdlambda} and repeating the arguments used in showing the estimate for $R_{g,\lambda}$ in Theorem \ref{Growing solutions Theorem}, then we get that there exists  $R_{d,\lambda}\in H^{1}(0,T;H^{-1}(M))\cap L^{2}(0,T;H^{1}(M))$  solving 
			\begin{align}\label{Equation for Rdlamdba after construction of pdi and Td}
				\begin{aligned}
						\mathcal{L}^{*}_{\varphi}\lb e^{i\psi}R_{d,\lambda}\rb(t,x)=-e^{i\psi}\mathcal{L}_{A,q}T_{d}(t,x),\ \ \ (t,x)\in M_{T}
				\end{aligned}
			\end{align}
			and $R_{d,\lambda}$ satisfies the following estimate \begin{align}\label{Estimate for Rdlambda}
				\lVert R_{d,\lambda}\rVert_{L^{2}(0,T;H^{1}_{\lambda}(M)}\leq C
			\end{align} for some constant $C>0$ independent of $\lambda$.
			Combining all these, we end up with proving the following theorem. 
			\begin{theorem}\label{Decaying solution Theorem}
				Let $M_{T}$, $\mathcal{L}_{A,q}$ and $\varphi$ be as before. Suppose $(D,g_{0})$ be a simple  manifold which is extension of $(M_{0},g_{0})$ in the sense that $M_{0}\subset D$ and $y_{0}\in  D$ is such that $(x_{1},y_{0})\notin M$ for all $x_{1}$. Now if $(r,\theta)$ denote the polar normal coordinates on $(D,g_{0})$, $(x_{1},r,\theta)$ denote the points in $M$ and $A_{1}$ and $A_{r}$ are components of $A$ in $x_{1}$ and $r$ coordinates respectively,  then for $\lambda$ large enough the following  equation  
				\begin{align*}
					\begin{aligned}
							\mathcal{L}^{*}_{A,q}v(t,x)=0,\ \ (t,x)\in M_{T}
					\end{aligned}
				\end{align*}
				has a solution taking the following form
				\begin{align}\label{eq:v in polar coordinates}
					\begin{aligned}
						v(t,x)=e^{-\lb \varphi-i\psi\rb(t,x)}\Big (T_{d}(t,x_{1},r,\theta)+R_{d,\lambda}(t,x_{1},r,\theta)\Big)
					\end{aligned}
				\end{align}
				where $\psi,T_{d}$ are given by \eqref{Expression for psi decaying case}, \eqref{Expression for Td decaying case} and $R_{d,\lambda}$ satisfies \eqref{Equation for Rdlamdba after construction of pdi and Td} and \eqref{Estimate for Rdlambda}. 
		\end{theorem}
		\section{Derivation of Integral identity and proof of Main Theorem}\label{Proof of main theorem}
		We use this section to derive an integral identity, which will be required to prove our main result. Later, using the geometric optics solutions constructed in Section \ref{Construction of GO solutions}, we conclude the proof of Theorem \ref{th:main theorem}. We  start by recalling
		\begin{align*}
			\mathcal{L}_{A,q} =  \partial_t - \sum_{j,k=1}^{n}\frac{1}{\sqrt{|g|}}\left(\partial_{x_j} +  A_j(t,x)\right)\left(g^{jk} \sqrt{|g|}(\partial_{x_{k}} + A_k(t,x))\right) + q(t,x)
		\end{align*}
		and \begin{align*}
			\mathcal{L}^*_{A, q} =   -  \partial_{t}-\sum_{j,k=1}^{n}\frac{1}{\sqrt{\lvert g\rvert}}(\partial_{x_{j}}-A_{j}(t,x))\lb\sqrt{\lvert g\rvert} g^{jk}(\partial_{x_{k}}- A_{k}(t,x))\rb+\overline{q}(t,x).
		\end{align*}
		For $l=1,2$, let $A^{(l)}$ and $q_l$ be as in Theorem \ref{th:main theorem}. Further assume that  $u_{l}$  is solution to the corresponding IBVP for $\mathcal{L}_{A^{(l)},q_{l}}$ given by \eqref{eq: Main equation} when $(A,q)=(A^{(l)},q_{l})$ for $l=1,2$, that is, for $l=1,2$, we have 
		\begin{align}\label{eq: equation for uj}
			\begin{aligned}
				\begin{cases}
					\mathcal{L}_{A^{(l)}, q_l}u_l(t, x)= 0, \ \  (t, x) \in M_T \\
					u_l(0, x)=  \phi(x), \ \ x \in M\\
					u_l(t, x)=  f(t, x), \ \ (t, x) \in \Sigma.
				\end{cases}
			\end{aligned}
		\end{align}
		Then $u:=  u_1 - u_2$, satisfies the following IBVP with zero initial and boundary conditions
		\begin{align}\label{eq: equation for the difference}
			\begin{aligned}
				\begin{cases}
					\mathcal{L}_{A^{(1)}, q_1}u(t, x) = \mathcal{Q} u_2(t, x) ,      \ \  (t, x) \in M_T \\
					u(0, x)=  0,   \ \ x \in M\\
					u(t, x) =  0, \ \  (t, x) \in \Sigma,
				\end{cases}
			\end{aligned}
		\end{align}
		where $\mathcal{Q} u_2(t, x):= \lb\lvert A^{(1)}\rvert_g^2 - \lvert A^{(2)}\rvert_g^2 \rb u_2 + 2  \left\langle A^{(1)} - A^{(2)}, \nabla_g u_2\right\rangle_{g} +  \delta_g \left(A^{(1)} - A^{(2)}\right)u_2 + (q_2 - q_1) u_2 $.
		To simplify the notation, let us denote by $\tilde{q}(t,x):=(\tilde{q}_{1}-\tilde{q}_{2})(t,x)$ and  $\tilde{A}(t,x):=(\tilde{A})_{1\leq j\leq n}:=(A^{(1)}-A^{(2)})(t,x)$ where  $\tilde{q}_i:= \lvert A^{(i)}\rvert^2_g + \delta_g A^{(i)} -q_i$,  for $i=1,2$, then with these notations  $\mathcal{Q} u_2$ becomes $$\mathcal{Q} u_2(t, x) =  2  \langle \tilde{A}(t,x), \nabla_g u_2(t,x)\rangle_{g} + \tilde{q}(t,x)u_2(t,x).$$ 
   Now since $\displaystyle \mathcal{Q} u_2\in L^2(M_T)$ therefore using Theorem 1.43 in \cite{Choulli_Book}
we have that there exists a unique solution  $\displaystyle u\in L^{2}(0,T;H^2(M))\cap H^1(0,T;L^2(M))$ to \eqref{eq: equation for the difference} with $\partial_\nu u\in L^2(0,T;H^{1/2}(\Sigma)).$
  Now if  $v(t, x)$ is a  solution to the adjoint operator of $\mathcal{L}_{A^{(1)}, q_1}$, given by 
\begin{align}\label{Adjoint equation integral identity}
		\begin{aligned}
					\mathcal{L}^*_{A^{(1)}, q_1}v(t, x)  = 0, \ \   (t, x) \in M_T,
			\end{aligned}
		\end{align}
then we observe that 
\begin{align*}
& \langle (\Lambda_{A^{(1)},q_1}-\Lambda_{A^{(2)},q_2})(\phi,f), v|_{\partial M^{*}_T}\rangle  = \langle \mathcal{N}_{A^{(1)},q_1}u_1 - \mathcal{N}_{A^{(2)},q_2}u_2, v|_{\partial M^{*}_T} \rangle \\
& \quad \quad = \int_{M_{T}}\left(-u_1\partial_t \overline{v}+\langle \nabla_g u_1,\nabla_g \overline{v}\rangle_g+2u_1\langle A^{(1)},\nabla_g \overline{v}\rangle_g +(\delta_g A^{(1)})u_1\overline{v}-\lvert A^{(1)}\rvert^2_g u_1 \overline{v}+q_1u_1\overline{v}  \right) \ dV_g dt \\
& \quad \quad \quad \quad-\int_{M}u_1(0,x)\bar{v}(0,x)dV_{g}\\
&\quad \quad \quad - \int_{M_{T}}\left(-u_2\partial_t \overline{v}+\langle \nabla_g u_2,\nabla_g \overline{v}\rangle_g+2u_2\langle A^{(2)},\nabla_g \overline{v}\rangle_g +(\delta_g A^{(2)})u_2\overline{v}-\lvert A^{(2)}\rvert^2_g u_2 \overline{v}+q_2u_2\overline{v} \right) \ dV_g dt\\
& \quad \quad \quad \quad  + \int_{M}u_2(0,x)\bar{v}(0,x)dV_{g}.
\end{align*}
Using integration by parts with $u|_{\Sigma} = 0, u|_{t=0} = 0$ and $v$ is solution to \eqref{Adjoint equation integral identity}, we get
\begin{align}{\label{eq: difference of Lambda}}
 \langle (\Lambda_{A^{(1)},q_1}-\Lambda_{A^{(2)},q_2})(\phi,f), v|_{\partial M^{*}_T}\rangle  = -  2  \int_{M_T}  \langle \tilde{A}(t,x), \nabla_g u_2(t,x)\rangle_{g} \bar{v}(t,x)\ dV_g d t - \int_{M_T} \tilde{q}(t,x)u_2(t,x) \bar{v}(t,x)\  dV_g dt.  
\end{align}

\noindent Multiplying equation \eqref{eq: equation for the difference} by  $\bar{v}(t, x)$ and integrate it over $M_T$, we get
		\begin{align*}
			\int_{M_T}  \mathcal{L}_{A^{(1)}, q_1}u(t, x) \bar{v}(t, x) \ dV_g dt=  2  \int_{M_T}  \langle \tilde{A}(t,x), \nabla_g u_2(t,x)\rangle_{g} \bar{v}(t,x)\ dV_g d t + \int_{M_T} \tilde{q}(t,x)u_2(t,x) \bar{v}(t,x)\  dV_g dt.
		\end{align*}
		Now use the integration by parts together with $u|_{\Sigma}= 0$, $u|_{t=0} = 0$, $\tilde{A}|_{\Sigma}=0$ 
		and the fact that  $v$ is a  solution to \eqref{Adjoint equation integral identity}, to  obtain the following  identity
		\begin{align}\label{eq: integral identity-I}
			\begin{aligned}
				& 2 \int_{M_T} \langle \tilde{A}(t,x), \nabla_g u_2(t,x)\rangle_{g} \bar{v}(t,x)\  dV_g dt + \int_{M_T} \tilde{q}(t,x)u_2(t,x) \bar{v}(t,x)\ dV_g  dt \\
				&\qquad  = -  \int_{\Sigma}  g^{jk}\nu_j \partial_{x_k} u(t,x) \bar{v}(t,x) \ dS_{g} dt + \int_M u(T,x) \bar{v}(T,x)\ dV_g. 
			\end{aligned}
		\end{align}
  From Equations \eqref{eq: difference of Lambda} and \eqref{eq: integral identity-I}, we have
  \begin{align*}
       \langle (\Lambda_{A^{(1)},q_1}-\Lambda_{A^{(2)},q_2})(\phi,f), v|_{\partial M^{*}_T}\rangle = \int_{\Sigma}  g^{jk}\nu_j \partial_{x_k} u(t,x) \bar{v}(t,x) \ dS_{g} dt - \int_M u(T,x) \bar{v}(T,x)\ dV_g. 
  \end{align*}
 Using \eqref{eq: Equal DN map }, we get $\partial_{\nu} u|_{\Sigma_{-,\epsilon/2}} = 0$ and $u|_{t=T} = 0$. Therefore, Equation \eqref{eq: integral identity-I} becomes
 \begin{align}\label{eq: integral identity}
			\begin{aligned}
				& 2 \int_{M_T} \langle \tilde{A}(t,x), \nabla_g u_2(t,x)\rangle_{g} \bar{v}(t,x)\  dV_g dt + \int_{M_T} \tilde{q}(t,x)u_2(t,x) \bar{v}(t,x)\ dV_g  dt \\
				&\qquad  = -  \int_{\Sigma \setminus \Sigma_{-,\epsilon/2}}  g^{jk}\nu_j \partial_{x_k} u(t,x) \bar{v}(t,x) \ dS_{g} dt. 
			\end{aligned}
		\end{align}
Let us define $J_1,J_2$ and $J_{3}$ by 
		\begin{align*}
			\begin{aligned}
				J_{1}:= 2 \int_{M_T} \langle \tilde{A}(t,x), &\nabla_g u_2(t,x)\rangle_{g} \bar{v}(t,x)\  dV_g dt,\ 
				J_{2}:=  \int_{M_T} \tilde{q}(t,x)u_2(t,x) \bar{v}(t,x)\ dV_g  dt \ \ 	\mbox{and}\\
				\ &J_{3}:= -  \int_{\Sigma\setminus \Sigma_{-,\epsilon/2}} \partial_{\nu}u(t,x) \bar{v}(t,x)\  dS_{g} dt. 
			\end{aligned}
		\end{align*}
		With these notations \eqref{eq: integral identity} becomes
		\begin{align}\label{Integral identity in compact form}
			J_{1}+J_{2}=J_{3}. 
		\end{align}
		Our next aim is to substitute the exponentially growing and decaying solutions constructed in section \ref{Construction of GO solutions},  for $u_{2}$ and $v$ respectively, in each term of equation \eqref{Integral identity in compact form}. Recall that $u_{2}$ satisfies 
		\begin{align*}
			\begin{aligned}
				\mathcal{L}_{A^{(2)},q_{2}}u_{2}=0,\ \mbox{in}\ M_{T}
			\end{aligned}
		\end{align*} and $v$ satisfies 
		\begin{align*}
			\mathcal{L}^{*}_{A^{(1)},q_{1}}v=0,\ \mbox{in}\ M_{T}
		\end{align*}
		therefore we choose the  expressions for solutions  $u_{2}$ and $v$ from  \eqref{eq:u in polar coordinates} and   \eqref{eq:v in polar coordinates} respectively, substitute in each term of  \eqref{Integral identity in compact form}. We start with the following calculations
		\begin{align*}
			\begin{aligned}
				&\Big\langle \tilde{A}(t,x),\nabla_{g}u_{2}(t,x)\Big\rangle_{g}\bar{v}=\Big(\Big\langle\tilde{A}(t,x), \nabla_g \left(\varphi+i\psi\right)\Big\rangle_{g}T_{g}(t,x)+ \Big\langle\tilde{A}(t,x), \nabla_g \left(\varphi+i\psi\right)\Big\rangle_{g}R_{g,\lambda}(t,x)\\
				&\qquad  + \Big\langle \tilde{A},\nabla_g T_{g}(t,x)\Big\rangle_{g}+\Big\langle \tilde{A},\nabla_g R_{g,\lambda}(t,x)\Big\rangle_{g}\Big)\Big(\Bar{T}_{d}(t,x)+\Bar{R}_{d,\lambda}(t,x)\Big).\\
				&=\Big\langle\tilde{A}(t,x), \nabla_g \left(\varphi+i\psi\right)\Big\rangle_{g}T_{g}(t,x)\Bar{T}_{d}(t,x)+\Big\langle\tilde{A}(t,x), \nabla_g \left(\varphi+i\psi\right)\Big\rangle_{g}T_{g}(t,x)\Bar{R}_{d,\lambda}(t,x)\\
				&\qquad +\Big\langle\tilde{A}(t,x), \nabla_g \left(\varphi+i\psi\right)\Big\rangle_{g}R_{g,\lambda}(t,x)\Bar{T}_{d}(t,x)+\Big\langle\tilde{A}(t,x)\nabla_g \left(\varphi+i\psi\right)\Big\rangle_{g}R_{g,\lambda}(t,x)\Bar{R}_{d,\lambda}(t,x)\\
				&\qquad \  +\Big\langle \tilde{A},\nabla_g T_{g}(t,x)\Big\rangle_{g}\Bar{T}_{d}(t,x)+\Big\langle \tilde{A},\nabla_g T_{g}(t,x)\Big\rangle_{g}\Bar{R}_{d,\lambda}(t,x)\\
				&\qquad \ +\Big\langle \tilde{A},\nabla_g R_{g,\lambda}(t,x)\Big\rangle_{g}\Bar{T}_{d}(t,x)+\Big\langle \tilde{A},\nabla_g R_{g,\lambda}(t,x)\Big\rangle_{g}\Bar{R}_{d,\lambda}(t,x)\\
				&:=\Big\langle\tilde{A}(t,x), \nabla_g \left(\varphi+i\psi\right)\Big\rangle_{g}T_{g}(t,x)\Bar{T}_{d}(t,x)+\mathcal{Z}_{1}(t,x). 
			\end{aligned}
		\end{align*}
		Similarly, we see that 
		\begin{align*}
			\begin{aligned}
				\tilde{q}(t,x)u_{2}(t,x)\bar{v}(t,x)=\tilde{q}(t,x)\big(\Bar{T}_{d}T_{g}(t,x)+\Bar{T}_{d}R_{g,\lambda}(t,x)+T_{g}\Bar{R}_{d,\lambda}(t,x)+\Bar{R}_{d,\lambda}(t,x)R_{g,\lambda}(t,x)\big):=\mathcal{Z}_{2}(t,x). 
			\end{aligned}
		\end{align*}
		Using the above expressions in definitions of $J_{1}$ and $J_{2}$, we get 
		\begin{align}\label{Sum of J1 and J2}
			\begin{aligned}
				J_{1}+J_{2}&=2\int_{M_{T}}\Big\langle\tilde{A}(t,x), \nabla_g \left(\varphi+i\psi\right)\Big\rangle_{g}T_{g}(t,x)\Bar{T}_{d}(t,x)\ dV_{g}dt\\
				&\qquad +2\int_{M_{T}}\mathcal{Z}_{1}(t,x)\ dV_{g}dt+\int_{M_{T}}\mathcal{Z}_{2}(t,x)\ dV_{g}dt. 
			\end{aligned}
		\end{align}
		Now using the expression for $v$ from \eqref{eq:v in polar coordinates} in the expression of $J_{3}$, we obtain
		\begin{align*}
			\begin{aligned}
				J_{3}&=-  \int_{\Sigma\setminus \Sigma_{-,\epsilon/2}} e^{-\lb \varphi+i\psi\rb}\partial_{\nu}u(t,x) \Bar{T}_{d}(t,x)\    dS_{g} dt 	- \int_{\Sigma\setminus \Sigma_{-,\epsilon/2}} e^{-\lb \varphi+i\psi\rb}\partial_{\nu}u(t,x) \Bar{R}_{d,\lambda}(t,x) \   dS_{g} dt.
			\end{aligned}
		\end{align*}
			We use the boundary Carleman estimate given in Theorem \ref{Boundary Carleman estimate Theorem} and follow  the arguments used in deriving Lemma $5.1$ in  \cite{Suman_Manmohan_IPI} to get the following estimate for  $J_{3}$
			\begin{align}\label{Estimate on J3}
				\lvert J_{3}\rvert\leq C\lambda^{1/2},\ \mbox{for some constant $C>0,$ independent of $\lambda$}. 
			\end{align}
			Using \eqref{Sum of J1 and J2} together with the estimate on $J_{3}$ given by \eqref{Estimate on J3} in \eqref{Integral identity in compact form}, we get 
			\begin{align*}
				\begin{aligned}
					\left \lvert 2\int_{M_{T}}\Big\langle\tilde{A}(t,x), \nabla_g \left(\varphi+i\psi\right)\Big\rangle_{g}T_{g}(t,x)\Bar{T}_{d}(t,x)\ dV_{g}dt\right\rvert
					&\leq  \left\lvert 2\int_{M_{T}}\mathcal{Z}_{1}(t,x)\ dV_{g}dt+\int_{M_{T}}\mathcal{Z}_{2}(t,x)\ dV_{g}dt\right\rvert +\lvert J_{3}\rvert\\
					&\leq C \left(\|\mathcal{Z}_{1}\|_{L^2(M_T)} + \|\mathcal{Z}_{2}\|_{L^2(M_T)} + |J_3|\right). 
				\end{aligned}
			\end{align*}
			Let $(x_{1},r,\theta)$ be the polar normal  coordinate on $(M,g)$ and $\tilde{A}_{1}$ and $\tilde{A}_{r}$ be components of $\tilde{A}$ in $x_{1}$ and $r$ coordinates respectively as in Theorem \ref{Growing solutions Theorem}. Then, the above estimate can be rewritten as 
			\begin{align}\label{Integral identity before estimates}
				\begin{aligned}
					\int_{\rr}\int_{\rr}\int_{S_{y}M_{0}}\int_{0}^{\tau_{+}(y_0,\theta)}\left(\tilde{A}_1+i\sqrt{1 - \beta^2}\tilde{A}_r\right) T_{g}(t,x_{1},r,\theta)\Bar{T}_{d}(t,x_{1},r,\theta)b(r, \theta)^{1/2} dx_1 drd\theta dt \\
					\leq \frac{C}{\lambda} \left(\|\mathcal{Z}_{1}\|_{L^2(M_T)} + \|\mathcal{Z}_{2}\|_{L^2(M_T)} + |J_3|\right)
				\end{aligned}
			\end{align}	
			where we used the fact $dV_g = b(r, \theta)^{1/2}dx_1 drd\theta$ in polar normal coordinates on $(M,g)$ and $\tau_{+}(y_0,\theta)$ is length of the geodesic in $M$, starting at $y_{0}$ in the direction of $\theta$.
			After using the estimates on $R_{d,\lambda}$ and $R_{g,\lambda}$ together with the expressions for $T_{d}$ and $T_{g}$  from Theorems \ref{Growing solutions Theorem} and \ref{Decaying solution Theorem}, we get that 
			\begin{align*}
				\begin{aligned}
					\lVert \mathcal{Z}_{i}\rVert_{L^{2}(M_{T})}\leq C, \ \mbox{for $i=1,2$ and constant $C>0$ independent of $\lambda$}. 
				\end{aligned}
			\end{align*}
			Using this estimate along with equation \eqref{Estimate on J3} in equation \eqref{Integral identity before estimates} and taking $\lambda\rightarrow \infty$, we get 
			\begin{align*}
				\int_{\rr}\int_{\rr}\int_{S_{y}M_{0}}\int_{0}^{\tau_{+}(y_0,\theta)}\left(\tilde{A}_1+i\sqrt{1 - \beta^2}\tilde{A}_r\right) \left(\phi(t)\right)^2e^{i\mu\left(\sqrt{1-\beta^{2}}\right)x_{1}}e^{-\mu r}e^{i\Phi(t, x_{1},r,\theta)}\left(h(\theta)\right)^2dx_1 drd\theta dt = 0,
			\end{align*}
			where  $\Phi(t, x_{1},r,\theta):= \left(\Phi_{1}-\overline{\Phi}_{2}\right)(t, x_{1},r,\theta)$ with $\Phi_1$ and $\Phi_2$ satisfying equations \eqref{Equation for Phi} and \eqref{Equation for Phi2 decaying case} respectively.
			As this relation is true for all cutoff functions $\phi \in C^\infty_c(0, T)$
			therefore we get
			\begin{align}\label{Final integral identity}
				\int_{\rr}\int_{S_{y}M_{0}}\int_{0}^{\tau_{+}(y_0,\theta)} \left(\tilde{A}_1+i\sqrt{1 - \beta^2}\tilde{A}_r\right)(t,x_{1},r,\theta)e^{i\mu\left(\sqrt{1-\beta^{2}}\right)x_{1}}e^{-\mu r}e^{i\Phi(t, x_{1},r,\theta)}\left(h(\theta)\right)^2   dx_1 dr d \theta= 0,
			\end{align}
			for all  $t \in (0, T)$ and $h\in C^{\infty}(S_{y_{0}}(D))$. Now following the arguments from \cite[see Section 6] {Ferreira_Kenig_Salo_Uhlmann_invention} we  obtain that there exists  a $\Psi \in W_0^{2,\infty}(M_T)$ such that $
			\tilde{A}(t,x) =  \nabla_g \Psi(t, x) \ \mbox{for}\  (t,x)\in M_{T}$. This proves the required uniqueness for the convection term.\\

			\noindent Next, we prove the uniqueness of density coefficient $q$. To prove this, We replace the pair $(A^{(1)},q_1)$ by $(A^{(3)},q_3)$, by taking $A^{(3)} = A^{(2)}$ in $M_{T}$, where $A^{(3)}(t,x) = A^{(1)}(t,x) - \nabla_g\Psi(t,x)$ and $q_3(t,x) = q_1(t,x) - \partial_t\Psi(t,x)$. From Proposition \ref{Gauge_equivalent_Prop} and Equation \eqref{eq: Equal DN map }, we get $\Lambda_{A^{(3)},q_3} = \Lambda_{A^{(2)},q_2}$. Using this in Equation \eqref{eq: integral identity}, 
			we get
			\begin{align*}
				\begin{aligned}
					\int_{M_T} (q_2 - q_3)(t,x)u_2(t,x) \bar{v}(t,x)\ dV_g  dt = -  \int_{\Sigma\setminus \Sigma_{-,\epsilon/2}}  g^{jk}\nu_j \partial_{x_k} u(t,x) \bar{v}(t,x)\  dS_{g} dt,
				\end{aligned}
			\end{align*}
   Again, we use the explicit expressions of $u_2$ and $v$ from Theorems \ref{Growing solutions Theorem} and \ref{Decaying solution Theorem} and take $\lambda\rightarrow \infty$ together with the estimate $\displaystyle \lVert \mathcal{Z}_{2}\rVert_{L^{2}(M_{T})}\leq C/\lambda,$ for some constant $C>0,$ independent of $\lambda$, to end up with getting 
			\begin{align*}
				\int_{\rr}\int_{\rr}\int_{S_{y}M_{0}}\int_{0}^{\tau_{+}(y_0,\theta)} q(t,x_{1},r,\theta)\left(\phi(t)\right)^2e^{i\mu\left(\sqrt{1-\beta^{2}}\right)x_{1}}e^{-\mu r}e^{i\Phi(t, x_{1},r,\theta)}\left(h(\theta)\right)^2 dt dx_1 dr d\theta  = 0,
			\end{align*}
			where $q(t,x_{1},r,\theta):=\lb q_{2}-q_{3}\rb (t,x_{1},r,\theta)$ is assumed to be zero outside $M_{T}$. 
			Finally by varying $\phi \in C^\infty_c(0, T)$, $h\in C^{\infty}(S_{y_{0}}(D))$ and taking $\Phi\equiv  0$ which is possible since $A^{(3)} = A^{(2)}$ in $M_{T}$  and  if $\Phi_1$ solves  \eqref{Equation for Phi} then we can choose $\Phi_2 =  \overline{\Phi}_1$ which solves \eqref{Equation for Phi2 decaying case}, we get that 
			\begin{align*}
				\int_{\rr}\int_{0}^{\tau_{+}(y_0,\theta)}q(t,x_{1},r,\theta)e^{i\mu\left(\sqrt{1-\beta^{2}}x_{1}+i r\right)}dx_1 dr= 0, \ \mbox{for all}\ \theta \in S_{y_{0}}(D), \ \mu>0,\ \beta\in \lb\frac{1}{\sqrt{3}},1\rb \   \mbox{and} \ t \in (0, T). 
			\end{align*}
			Now following \cite{KIna_Oksanen_IMRN}, we  vary $y_{0}\in D$ such that $(x_{1},y_{0})\notin M$ for all $x_{1}$ to get 	 $q\equiv 0,$ in $M_{T}$  which gives  $q_{1}(t,x) - q_{2}(t,x) = \partial_t\Psi(t,x)$ for $(t,x)\in M_{T}$. 	This completes the proof of the main theorem.

\section*{Acknowledgments}
RM was partially supported by SERB SRG grant No. SRG/2022/000947. MV is supported by ISIRD project from  IIT Ropar and Start-up Research Grant SRG/2021/001432  from the Science and Engineering Research Board, Government of India. AP is supported by UGC, Government of India, with a research fellowship.

		\end{document}